\documentclass[onefignum,onetabnum]{siamart171218}
\usepackage{bbm}
\usepackage{subfig}
\usepackage{amsfonts}
\usepackage{graphicx}
\usepackage{epstopdf}
\usepackage{algorithmic}
\DeclareMathOperator*{\argmin}{arg\,min}

\ifpdf
  \DeclareGraphicsExtensions{.eps,.pdf,.png,.jpg}
\else
  \DeclareGraphicsExtensions{.eps}
\fi

\newsiamremark{remark}{Remark}
\newsiamremark{hypothesis}{Hypothesis}
\crefname{hypothesis}{Hypothesis}{Hypotheses}
\newsiamthm{claim}{Claim}
\crefname{claim}{Claim}{Claim}
\newcommand{\alexander}[1]{\textcolor{blue}{ #1}}

\newcommand{\evo}[1]{\mathcal{S}_{#1}}

\headers{On The VIIM}{On The VIIM}

\title{On the Voronoi Implicit Interface Method\thanks{\textbf{Funding:} Alexander Zaitzeff and Selim Esedoglu gratefully acknowledge support from the NSF grant DMS-1719727. Krishnna Garikipati acknowledges NSF grants DMREF-1436154 and DMREF-1729166.}}

\author{Alexander Zaitzeff\thanks{Department of Mathematics, University of Michigan, Ann Arbor, MI 48109, USA
  (\email{azaitzef@umich.edu}, \email{esedoglu@umich.edu}).}
\and Selim Esedoglu\footnotemark[2]
\and Krishna Garikipati\thanks{Departments of Mechanical Engineering, and Mathematics, Michigan Institute for Computational Discovery \& Engineering, University of Michigan, Ann Arbor, MI 48109, USA (\email{krishna@umich.edu}).} }

\usepackage{amsopn}

\ifpdf
\hypersetup{
  pdftitle={On The VIIM},
  pdfauthor={Zaitzeff, Esedoglu, Garikipati}
}
\fi

\begin{document}

\maketitle

\begin{abstract}
We present careful numerical convergence studies, using parameterized curves to reach very high resolutions in two dimensions, of a level set method for multiphase curvature motion known as the Voronoi implicit interface method.
Our tests demonstrate that in the unequal, additive surface tension case, the Voronoi implicit interface method does not converge to the desired limit.
We then present a variant that maintains the spirit of the original algorithm, and appears to fix the non-convergence.
As a bonus, the new variant extends the Voronoi implicit interface method to unequal mobilities.
\end{abstract}



\section{Introduction}

The Voronoi implicit interface method (VIIM)~\cite{saye2012analysis} is a type of level set method \cite{osher1988fronts} that is particularly suited to the treatment of multiphase flows. It has been demonstrated on problems that incorporate surface tension, such as dynamics of bubble clusters, and the motion of grain boundaries. Often, the hardest aspect of designing numerical schemes for such problems is ensuring that the correct angle conditions at triple junctions -- in terms of the surface tensions of the interfaces meeting there -- are satisfied. The purpose of this paper is to investigate whether the VIIM in fact attains the correct angle conditions at junctions. To that end, and to study this essential difficulty in isolation, we focus on problems where surface tension is the sole driving force: multiphase motion by mean curvature.
This evolution arises in many applications, including grain boundary motion in polycrystalline materials during annealing \cite{mullins1956two} and image segmentation algorithms in computer vision \cite{mumford1989optimal}.
Once the correct behavior at junctions has been verified, additional driving forces, e.g. from bulk effects, fluid flow, etc., can be incorporated in the VIIM and similar algorithms with ease.\\

The first contribition of this paper is to present convergence tests at very high resolutions for the VIIM. To our knowledge, this is the first time that such highly accurate convergence tests have been done for this method. One of our main results is that the VIIM does not, in general, converge to the correct solution when interfaces with unequal surface tensions meet at a triple junction. Our second contribution in this paper is to present a modification of the Voronoi implicit interface method that fixes the non-convergence and ensures that the correct angle conditions are attained at all triple junctions.

The paper is organized as follows:
\begin{itemize}
\item In \cref{sec:cm} we briefly review multiphase motion by mean curvature and recall its variational formulation.
\item In \cref{sec:viim} we recall the Voronoi implicit interface method.
\item In \cref{sec:main} we present our implementation of the VIIM using parameterized curves.
We then present high resolution convergence tests indicating that the VIIM does not converge to the correct solution when the surface tensions are unequal.
\item In \cref{sec:thres} we present convergence tests using the parametrized curve representation of another recent algorithm, known as threshold dynamics, the convergence of which is supported by many recent results in the literature. This verifies the validity and accuracy of our parametrized curve implementation.
\item In \cref{sec:alg} we present a modification to the VIIM, which may be called the dictionary mapping implicit interface method (DMIIM), that does converge in the unequal surface tension case.
After discussing its relation to the VIIM, we subject this new variant to the same careful numerical convergence studies, again via an implementation on parametrized curves  to reach very high resolutions, and thus demonstrate its accuracy and convergence.
We then present further examples of the new algorithm on implicitly defined interfaces on a grid.
\end{itemize}
The code for sections 4 through 6 is publicly available, and can be found at \url{https://github.com/AZaitzeff/DMIIM}.
\section{Multiphase Motion by Mean Curvature}
\label{sec:cm}
In this paper, we will be concerned exclusively with gradient descent dynamics for energies of the form
\begin{align}
\label{energy}
E(\Sigma_1,\ldots,\Sigma_n)=\sum_{i\neq j} \sigma_{ij} \text{Area}(\Gamma_{ij}).
\end{align}
where $\Gamma_{ij} = (\partial\Sigma_i) \cap (\partial\Sigma_j)$ are the interfaces between the {\it phases} $\Sigma_1,\ldots,\Sigma_n$ that partition a domain $D\subset\mathbb{R}^d$, $d \geq 2$:
$$ \Sigma_{i} \cap \Sigma_{j} = (\partial\Sigma_i) \cap (\partial\Sigma_j) \mbox{ for any } i\not= j \mbox{, and } \bigcup _{i=1}^N \Sigma_i = D.$$
Note that if two grains are not neighbors, the intersection is empty and the corresponding term in the sum of \cref{energy} drops out. The positive constants $\sigma_{ij}=\sigma_{ji}$ are known as {\it surface tensions} (or surface energy density).
They need to satisfy the triangle inequality
\[
\sigma_{ij}+\sigma_{ik}\geq \sigma_{jk} \text{  for any distinct $i,j$ and $k$}
\] 
for well-posedness of the model (\ref{energy}).
Multiphase mean curvature motion arises as $L^2$ gradient descent on energies of this form, and can be described as follows:
\begin{enumerate}
\item At any point $p \in \Gamma_{ij}$ away from triple junctions where the interface is smooth, the normal speed, denoted $v_{\perp}(p)$, is given by 
\begin{align}
\label{nvel}
v_{\perp}(p)=\mu_{ij}\sigma_{ij}\kappa_{ij}(p).
\end{align}
Here, $\kappa_{ij}(p)=\kappa_{ji}(p)$ is the {\it mean curvature} at point $p$. The positive constants $\mu_{ij}$ are called {\it mobilities}.
Unless otherwise stated, we will take $\mu_{ij}=1$. 
\item Let a {\it triple junction} be formed by the meeting of three phases $\Sigma_1$, $\Sigma_2$, and $\Sigma_3$. They are points in two dimensions and occur along curves in three dimensions. Let $\theta_i$ be the angle between $\Gamma_{ij}$ and $\Gamma_{ik}$ at the junction.
Then:
\begin{align}
\label{herring}
\frac{\sin{\theta_1}}{\sigma_{23}}=\frac{\sin{\theta_2}}{\sigma_{13}}=\frac{\sin{\theta_3}}{\sigma_{12}}
\end{align}
has to hold.
This is known as the Herring angle condition~\cite{herring1999surface}. 
\end{enumerate}

Until now, the VIIM has been described only for the very special class of surface tensions known as {\it additive} surface tensions:
We call $\sigma_{ij}$ additive if they can be written as
\begin{equation*}
    \sigma_{ij} = \frac{1}{2}(\sigma_i + \sigma_j)
\end{equation*}
for some constants $\sigma_1,\sigma_2,\ldots,\sigma_n \geq 0$.
Additive surface tensions thus have $n$ degrees of freedom, constituting therefore a very small subclass of physically relevant surface tensions, which require $\binom{n}{2}$ degrees of freedom to fully specify.


\section{The Voronoi Implicit Interface Method}
\label{sec:viim}

\begin{figure}
  \begin{center}
\includegraphics[width=.3\textwidth]{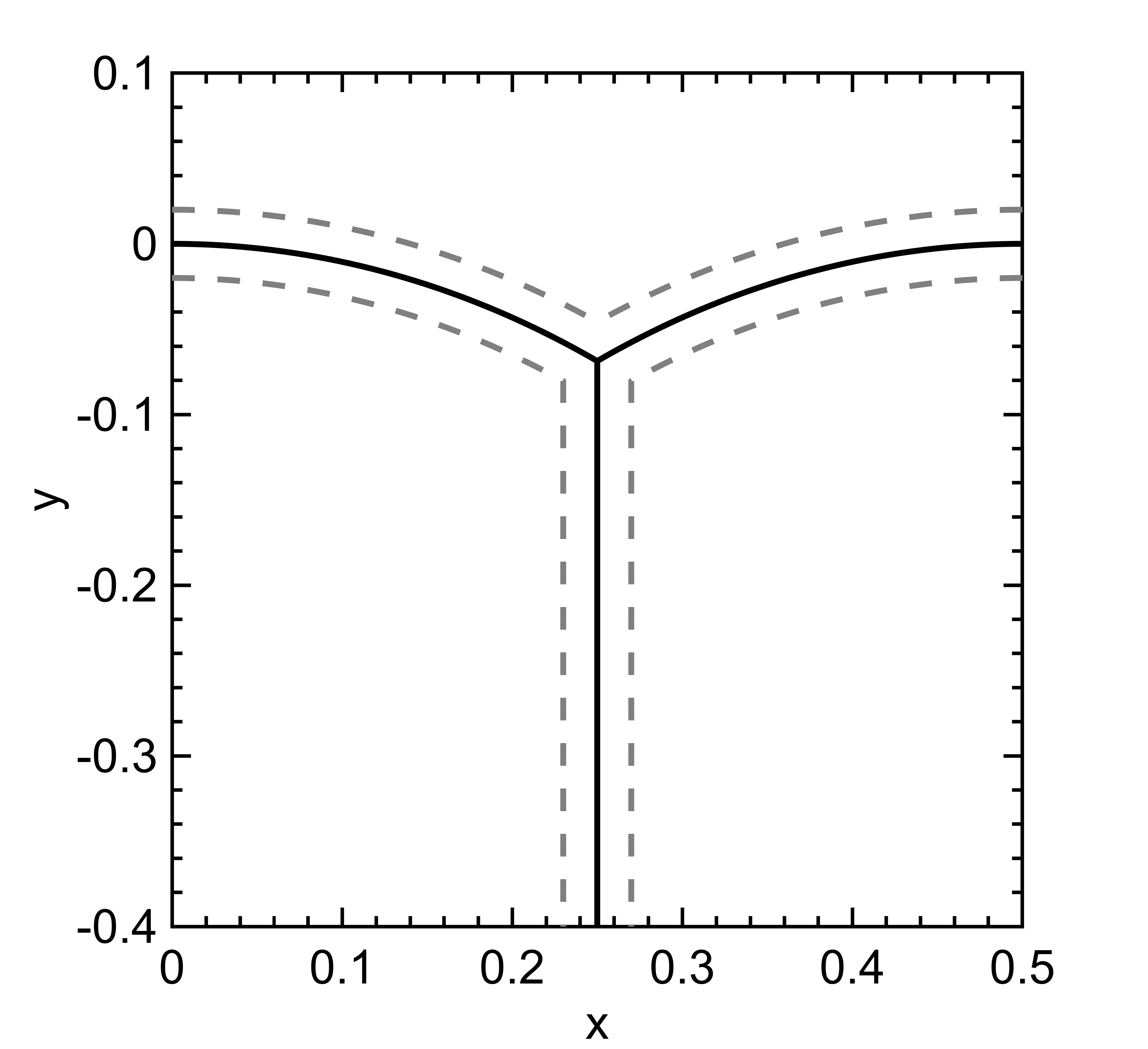} %
\includegraphics[width=.3\textwidth]{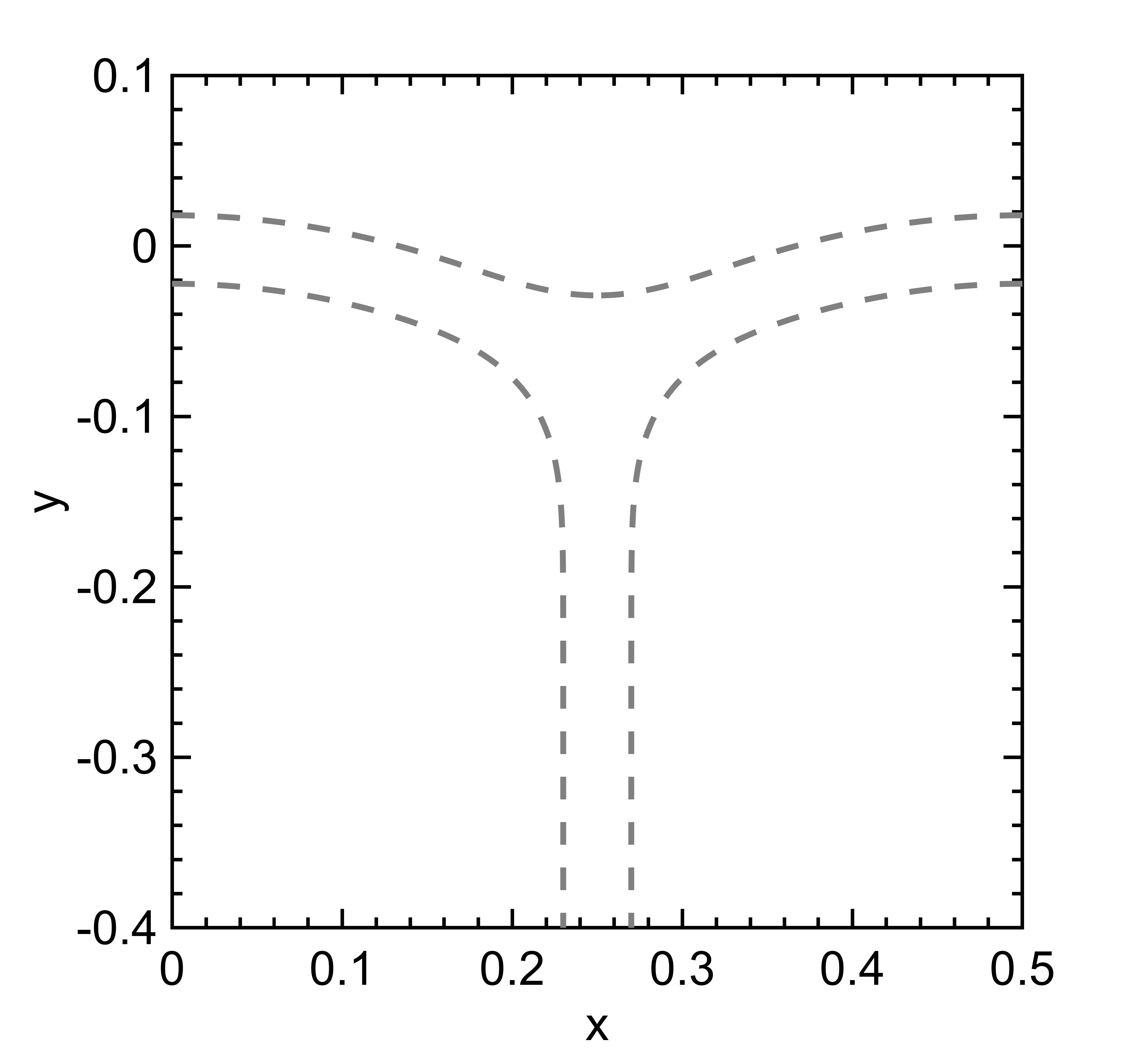}%
\includegraphics[width=.3\textwidth]{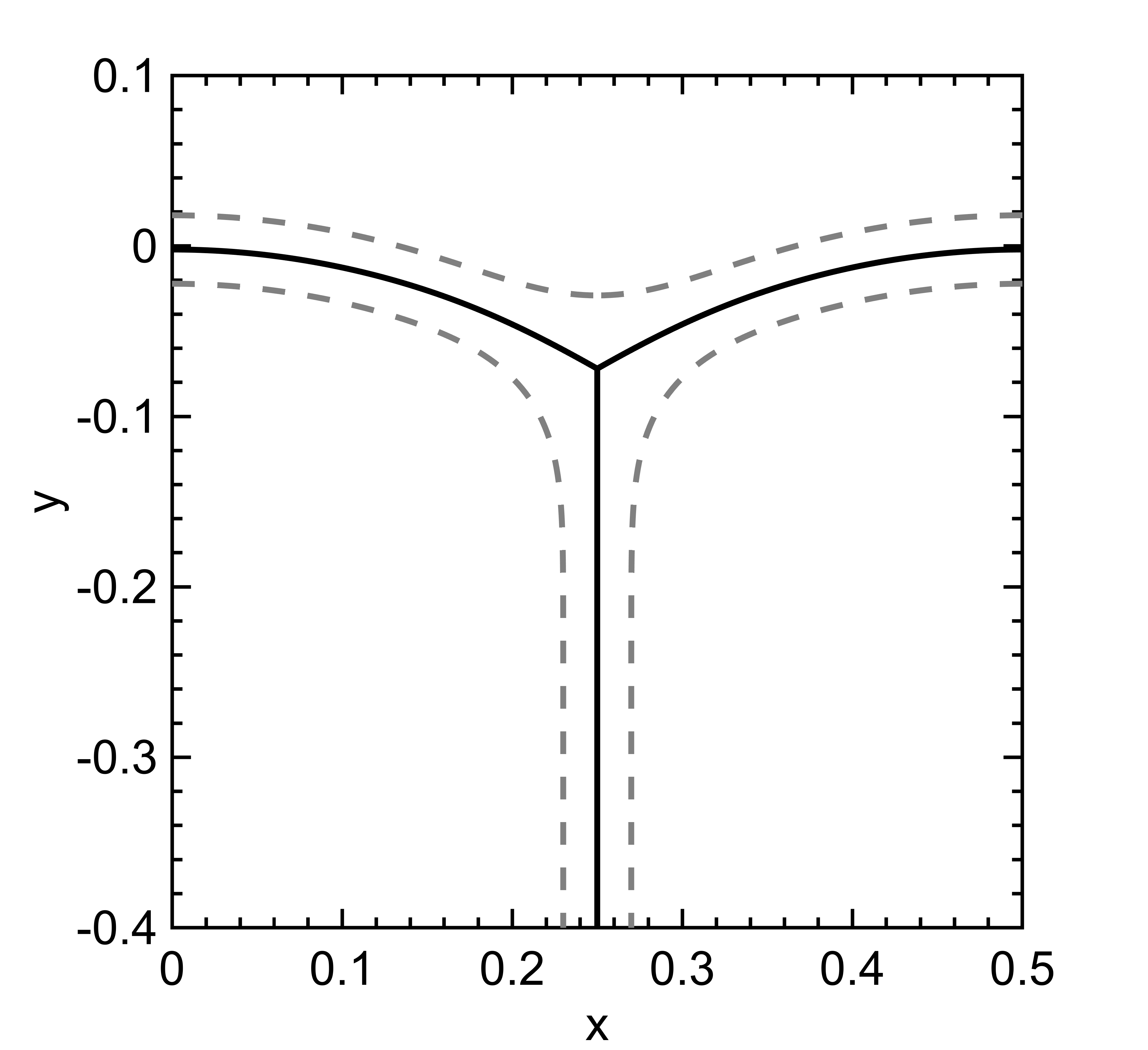}%
    \end{center}
    \caption{\footnotesize How the VIIM works: On the left with have the $\epsilon$ level sets (dotted lines), in the center we have the sets after being evolved by \ref{eq:ls}, the figure on the right shows the new interfaces after Voronoi reconstruction (solid line).}
    \label{fig:viim}
\end{figure}
We recall the level set formulation \cite{osher1988fronts} of motion by mean curvature that gives the current configuration, $\partial\Sigma(t)$, of the boundary with initial configuration $\partial\Sigma(0)$ of a given initial set $\Sigma(0)\subset D$ in the two phase setting.
Let $\phi:D\times \mathbb{R}^+ \to \mathbb{R}$ be a level set function for $\Sigma(0)$ so that $\phi(x,0) > 0$ for $x\in\Sigma^\circ(0)$ and $\phi(x,0)<0$ for $x\in\Sigma^c(0)$.
Often, a particularly convenient choice of level set function to represent the boundary $\partial\Sigma$ of a set $\Sigma$ is the {\it signed distance function},
\begin{align}
\label{eq:distance}
d_{\Sigma}(\mathbf{x})= \begin{cases}
   \displaystyle \min_{\mathbf{z} \in \partial \Sigma} ||\mathbf{x}-\mathbf{z}||_2 & \mathbf{x}\in \Sigma \\
  -\displaystyle \min_{\mathbf{z} \in \partial \Sigma} ||\mathbf{x}-\mathbf{z}||_2 & \mathbf{x} \not\in \Sigma. 
  \end{cases}
\end{align}
The process of ``reinitializing'' a level set function $\phi$ by replacing it with the signed distance function to its $0$-level set, $d_{ \{ x \, : \, \phi(x) \geq 0 \} }(x)$, is known as ``redistancing'' in the level set literature, and is a common operation, typically applied only sporadically to prevent $\phi$ from becoming too steep or too flat (see, for example,~\cite{sethian1999level}).
In any case, if $\phi(x,t)$ solves the well-known PDE
\begin{align}
\label{eq:ls}
\phi_t-\nabla\cdot\bigg(\frac{\nabla \phi}{|\nabla \phi|}\bigg)|\nabla \phi|=0.
\end{align}
and we set $\Sigma(t) = \{ \mathbf{x} \in D \, : \, \phi(\mathbf{x},t) \geq 0 \}$, then $\partial\Sigma(t)$ evolves by motion by mean curvature.\\ 

There have been multiple algorithms proposed in the literature to extend the level set formulation of mean curvature motion to the multiphase setting. We note, in particular, the level set method of \cite{chopp}, the variational level set method of \cite{zhao1996variational}, and the distance function based diffusion generated motion of \cite{elsey2009diffusion}, \cite{esedog2010diffusion}, and \cite{elsey2013}.
The latter three contributions alternate diffusion by the linear heat equation applied individually to the level set functions of the phases (so that they are decoupled at this stage), a simple pointwise redistribution step that couples the phases, and redistancing on the new level set functions, to generate the desired multiphase evolution; in this sense, they are a cross between the convolution generated motion of \cite{merriman1994motion} and the level set method.\\

In \cite{sspnas}, Saye and Sethian introduced a variant of the algorithm in \cite{elsey2009diffusion}.
This new version also alternates redistancing and decoupled evolution of the level sets of individual phases with a pointwise redistribution step that imposes the requisite coupling.
The key differences are: {\bf 1.} the decoupled motion of the level sets is by the nonlinear PDE \cref{eq:ls} vs.\ the linear heat equation, and {\bf 2.} the redistribution step takes place after (vs.\ before) redistancing of the individually evolved level set functions.
An additional novelty, mostly for convenience, is an innovative step to enable evolution of all the level set functions concurrently, by applying \cref{eq:ls} to the unsigned distance function of the union of the $\epsilon>0$ super-level sets of the phases, 
\[\varphi_\epsilon(\mathbf{x})=    \displaystyle \min_{\mathbf{z} \in \cup_i \partial \Sigma_i} ||\mathbf{x}-\mathbf{z}||_2-\epsilon.\]

Although only the equal surface tension case ($\sigma_{ij} = 1$ for all $i\not= j$) of multiphase motion by mean curvature was considered in the original paper \cite{sspnas}, in a subsequent contribution \cite{saye2012analysis}, Saye and Sethian proposed an extension of their method to certain (additive) unequal surface tension networks.
This is a very small subset of all surface tensions allowed by model (\ref{energy}).
Moreover, the extension in \cite{saye2012analysis} takes all mobilities to be equal, again a vast restriction over (\ref{nvel}).
One of the original motivations for the present study was to investigate if the algorithm could be extended to the full generality of model (\ref{energy}) \& (\ref{nvel}).
\\

We will introduce some notation to represent various steps of the VIIM as described in \cite{sspnas}, including the extension to additive surface tensions in \cite{saye2012analysis}.
To that end, first define the function $\mathcal{S}_{\Delta t}$ by
\begin{equation*}
(d_{\Sigma_1(\sigma_1 \Delta t,\epsilon)},d_{\Sigma_2(\sigma_2 \Delta t,\epsilon)},\ldots,d_{\Sigma_n(\sigma_n \Delta t,\epsilon)})=\evo{\Delta t}(\Sigma_1,\Sigma_2,\ldots,\Sigma_n)
\end{equation*}
where $\Sigma_i(t,\epsilon) = \{ \mathbf{x} \, : \, \phi_i(\mathbf{x},t) \geq \epsilon \}$ denotes the $\epsilon$-super level set of the solution $\phi_i(\mathbf{x},t)$ of mean curvature flow equation (\ref{eq:ls}) at time $t$, starting from the initial condition $\phi_i(\mathbf{x},0) = d_{\Sigma_i}(\mathbf{x})$. Here $\sigma_i$ is the surface tension associated with phase $\Sigma_i$ whereas the surface tension corresponding to the interface $\Gamma_{ij}$ is $\frac{1}{2}(\sigma_i+\sigma_j)$.\\

Next, the Voronoi reconstruction step of the VIIM (that reallocates points among the phases) will be represented by the function $R_v$, which maps an $n$-tuple of functions $(\phi_1,\ldots,\phi_n)$ to an $n$-tuple of sets $\Omega_1,\Omega_2,\ldots,\Omega_n$:
\begin{equation*}
(\Omega_1,\Omega_2,\ldots,\Omega_n) = R_v ( \phi_1, \phi_2, \ldots, \phi_n)
\end{equation*}
where
\begin{equation*}
\Omega_i = \Big\{ \mathbf{x} \, : \, \phi_i(\mathbf{x}) = \max_{j} \phi_j(\mathbf{x}) \Big\}.
\end{equation*}

With this notation, the evolution of a multiphase system by the VIIM at the $N$-th time step with time step size $\Delta t$, $T=N\Delta t$,  is given by
\begin{equation}
\label{eq:VIIM}
(R_v \circ \evo{\Delta t})^N (\Sigma_1,\Sigma_2,\ldots,\Sigma_n).
\end{equation}
In \cite{sspnas,saye2012analysis}, no extension of the algorithm is given to the far more general case of $\binom{n}{2}$ surface tensions; moreover, the mobilities of all the interfaces are assumed to be $1$, with again no indication given for greater generality.
The method is summarized by \cref{alg:VIIMe} and illustrated in \cref{fig:viim}.\\
\begin{algorithm}
  \caption{The Voronoi Implicit Interface Method
    \label{alg:VIIMe}}
  \begin{algorithmic}[1]
  	\STATE Given $\Sigma^0_1,\Sigma^0_2,\ldots,\Sigma^0_n$.
    \STATE Let $N=T/\Delta t$.
      \FOR{$k \gets 1 \textrm{ to } N$}
      
          \STATE Evolve each $\phi_i(\cdot,0)=d_{\Sigma^{k-1}_i}$ by time $\sigma_i \Delta t$ by \ref{eq:ls} to obtain $\phi_i(\cdot,\sigma_i \Delta t)$.
          \STATE Build the signed distance functions $d_{\Sigma_i^{k-1}(\sigma_i \Delta t,\epsilon)}=d_{\{\mathbf{x}:\phi_i(\mathbf{x},\sigma_i \Delta t) \geq\epsilon\}}$. 
          \STATE Construct the new phases $\Sigma^k_i=\{\mathbf{x}:d_{\Sigma_i^{k-1}(\sigma_i \Delta t,\epsilon)}(\mathbf{x}) = \max_{j} d_{\Sigma_j^{k-1}(\sigma_j \Delta t,\epsilon)}(\mathbf{x})\}$
      \ENDFOR
  \end{algorithmic}
\end{algorithm}

Saye and Sethian state in \cite{saye2012analysis} that convergence to the desired motion is obtained by taking the double limit, $\displaystyle \lim_{\epsilon\to 0} \lim_{\Delta t \to 0}$ in (\ref{eq:VIIM}), with $N = \frac{T}{\Delta t}$.
However, they also discuss two other limiting procedures: the ``coupled'' limit, $\displaystyle \lim_{\epsilon=c\sqrt{\Delta t} \to 0^+}$, and the interchanged double limit, $\displaystyle \lim_{\Delta t \to 0} \lim_{\epsilon \to 0}$.
In the case of equal surface tensions, the authors present numerical convergence studies for each of these three limits.
In the case of unequal surface tension, they only cite qualitative evidence and only for the coupled limit.
In the next section, we present highly accurate, exhaustive numerical convergence studies of the VIIM, in the equal and unequal surface tension cases, for {\it all} of these limits.

\section{Testing the VIIM using Parameterized Curves}
\label{sec:main}
To carefully assess the convergence of the VIIM, we will implement it in $\mathbb{R}^2$ using {\it parametrized curves} to represent the interfaces $\Gamma_{ij} = (\partial\Sigma_i) \cap (\partial\Sigma_j)$, and test it on exact solutions away from topological changes.
This will allow us to reach resolutions not easily attainable with the practical implementation of the algorithm for arbitrary initial data via implicit (level set) representation of the interfaces on a uniform grid.
We stress that we are not advocating parametric representation in the context of the VIIM as a general numerical method, as it defeats the original purpose -- painless handling of topological changes -- of a level set based algorithm; we use it only to carry out reliable convergence studies.\\

Below, Section 4.1 recalls the well-known ``Grim Reaper'' exact solution of three phase motion by mean curvature, with equal and unequal surface tensions.
Section 4.2 discusses in detail how each step of the VIIM is implemented via parametrized curves (front tracking).
Finally, Section 4.3 presents the results of our numerical convergence study.

\subsection{``Grim Reaper'' Solution}
\begin{figure}
  
  \begin{center}
  \subfloat{\includegraphics[width=.25\textwidth]{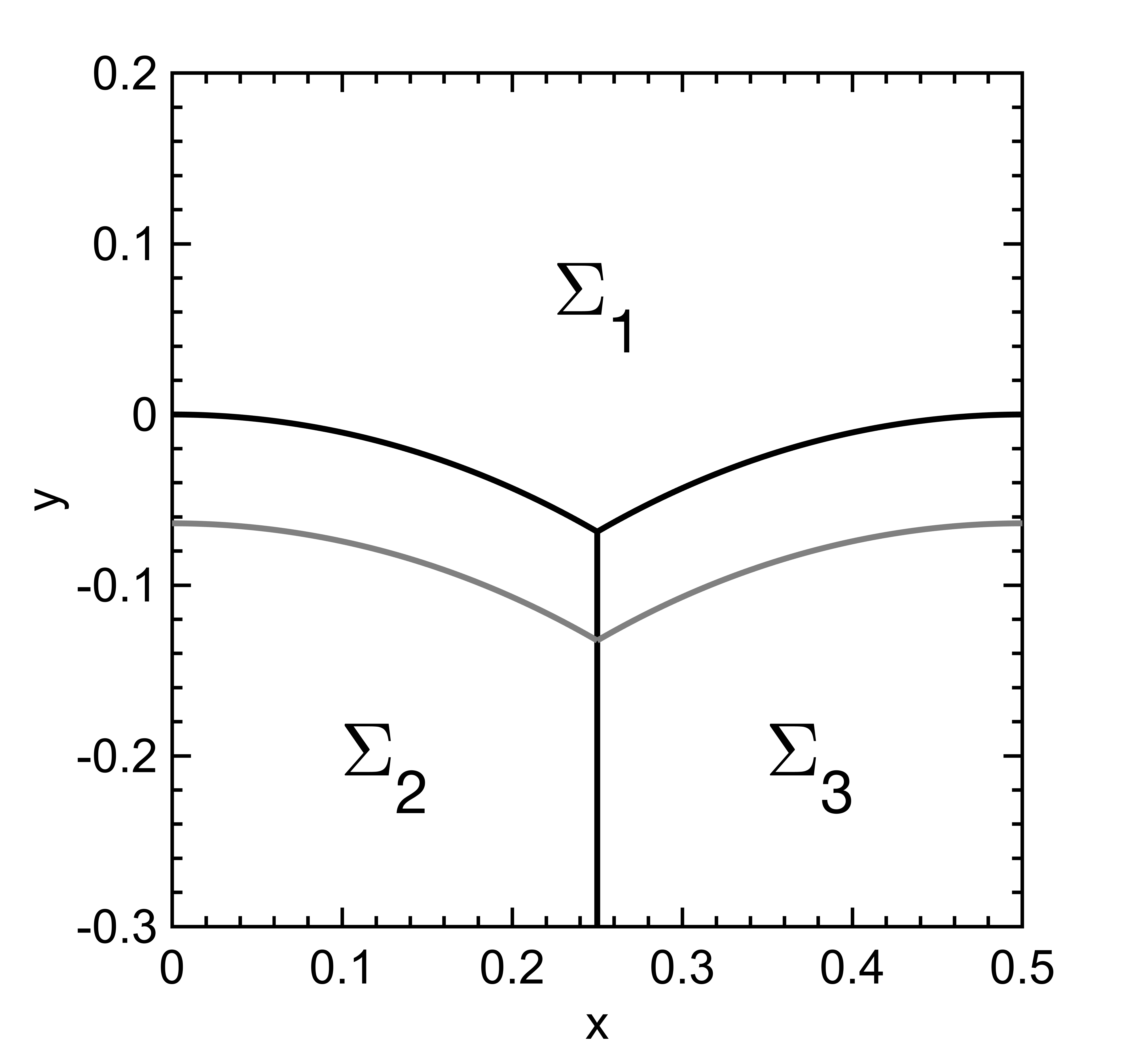}}
  \subfloat{\includegraphics[width=.25\textwidth]{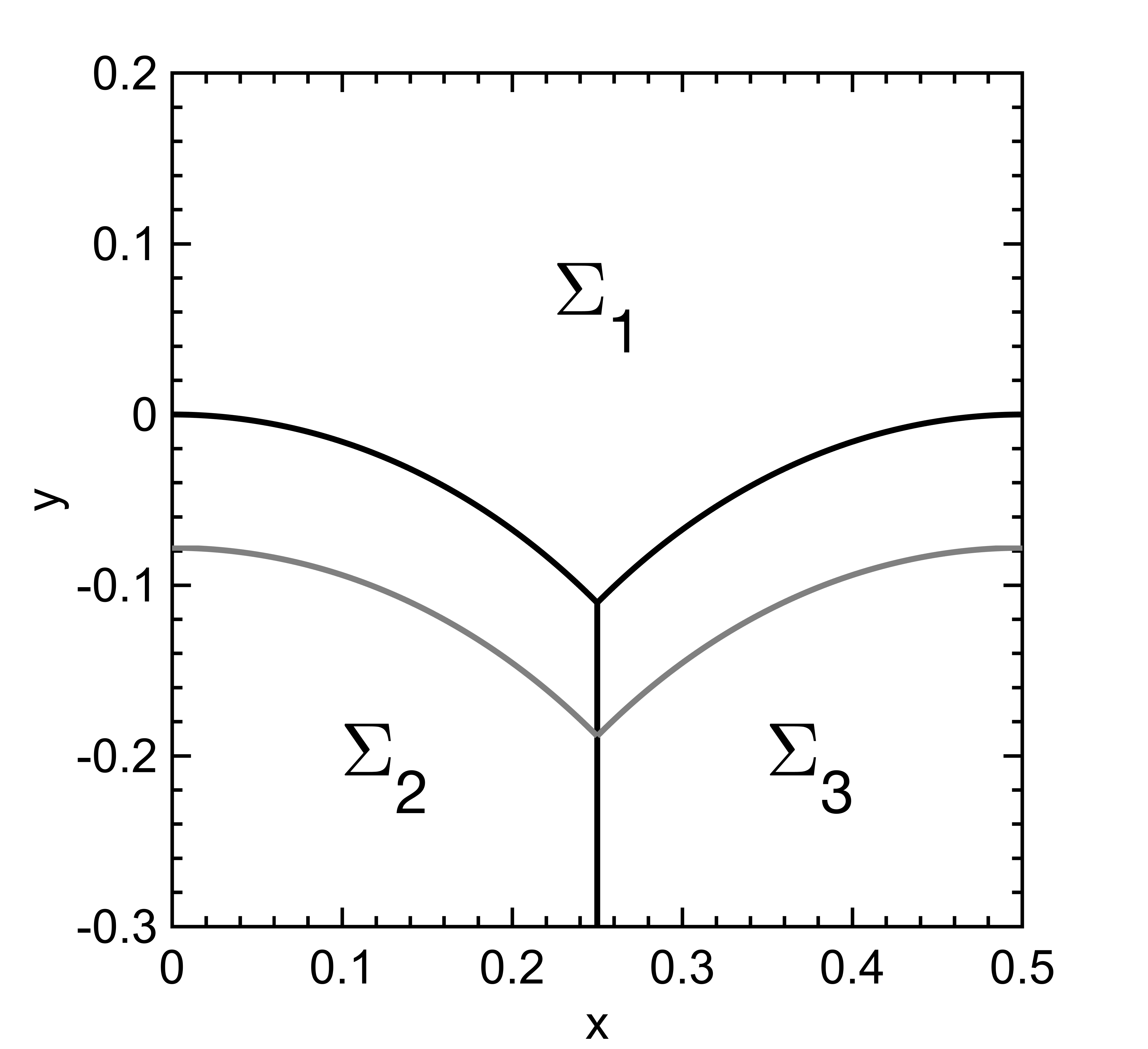}}
  \subfloat{\includegraphics[width=.25\textwidth]{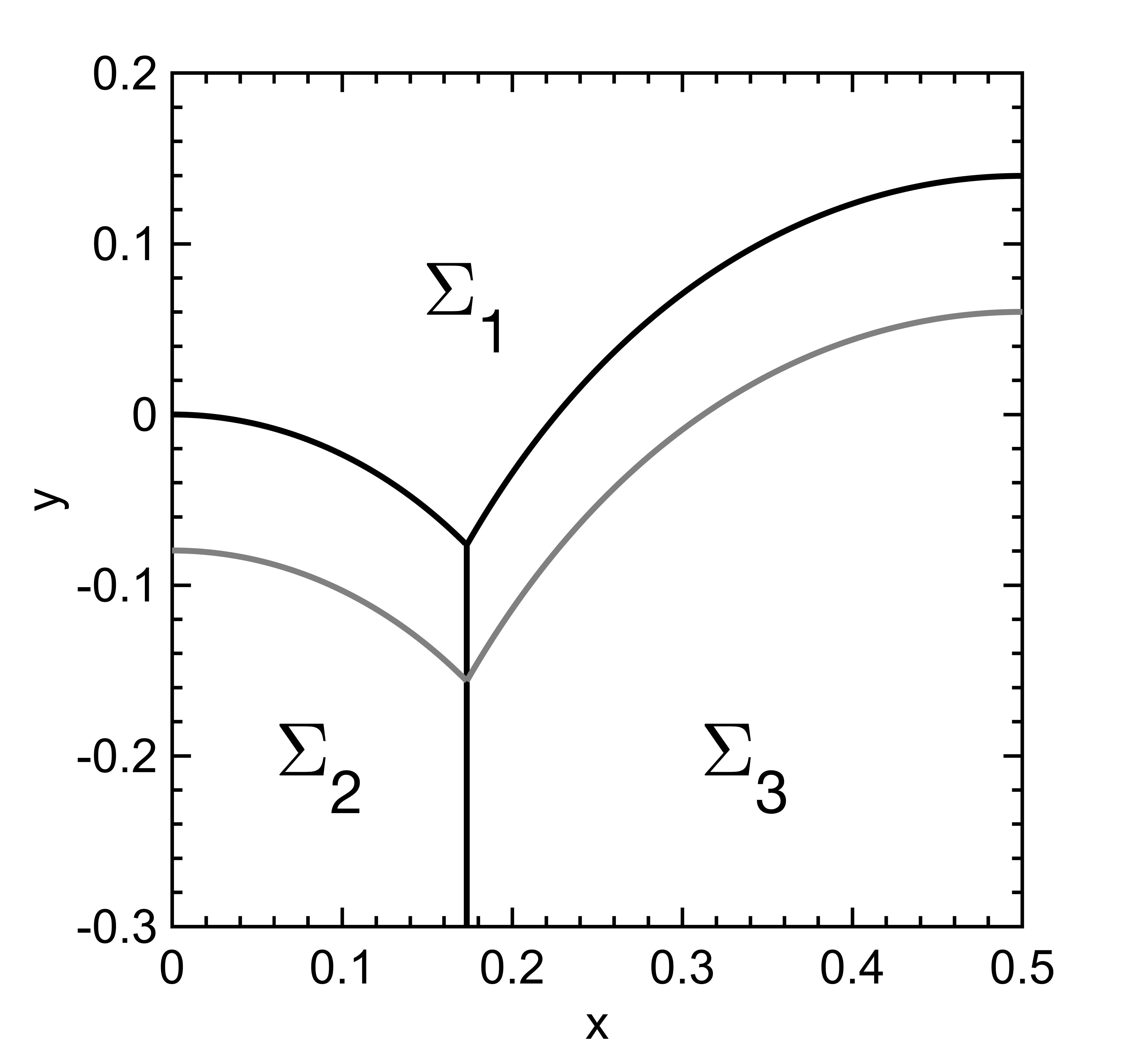}}
  \subfloat{\includegraphics[width=.25\textwidth]{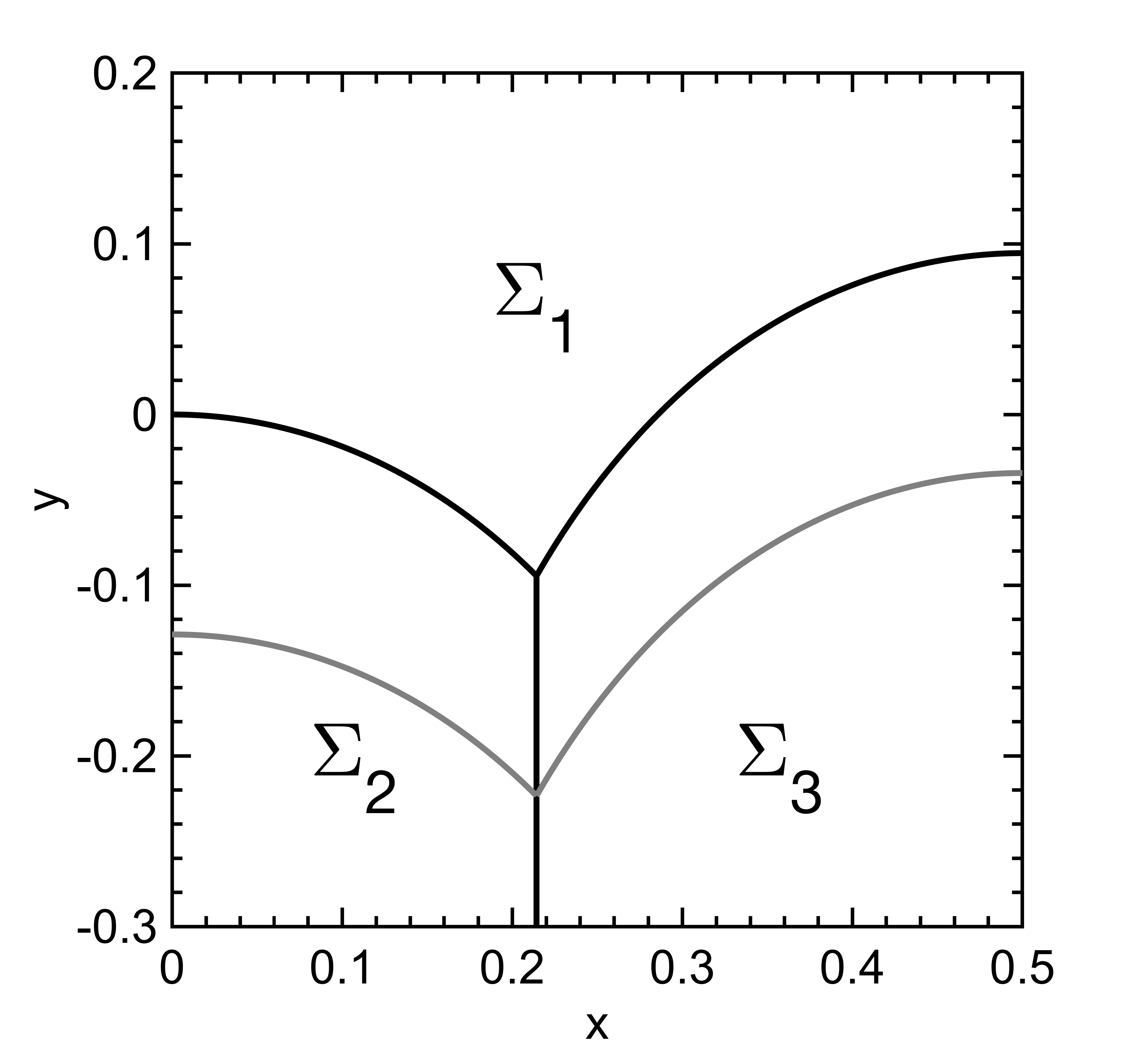}}
  
  \end{center}
  \caption{\footnotesize ``Grim Reaper'' exact solutions for angles (left to right) $(120^\circ,120^\circ,120^\circ)$, $(90^\circ,135^\circ,135^\circ)$,  $(75^\circ,135^\circ,150^\circ)$ with $\mu_i=1$ and  $(75^\circ,135^\circ,150^\circ)$  with $\mu_i=\frac{1}{\sigma_i}$. The black line is $t=0$ and the gray line is $t=\frac{18}{512}$. }
  \label{fig:gr75}
\end{figure}

``Grim Reapers'' are a family of exact solutions to three-phase motion by mean curvature that include unequal surface tension and mobility cases.
Our domain $D$ will be $\big[0,\frac{1}{2}\big] \times [-\frac{1}{2},1]\subset \mathbb{R}^2$, and we will impose Neumann boundary conditions: interfaces intersect $\partial D$ at right angles.
In all our examples, the interface $\partial\Sigma_1(t)$ will be given as the graph of a function $f=f(x,t)$, at least on $t\in[0,\frac{18}{512}]$ where we choose $18/512$ to allow the curves to travel a appreciable distance in our domain.
The three phases $\Sigma_1(t)$, $\Sigma_2(t)$, and $\Sigma_3(t)$ will be described in terms of $f(x,t)$ as follows: 
\begin{align*}
\Sigma_1(t)=&\bigg\{(x,y):  y \geq f(x,t)  \bigg\}\\
\Sigma_2(t)=&\bigg\{(x,y): x\leq \beta \text{ and } y \leq f(x,t)\bigg\}\\
\Sigma_3(t)=&\bigg\{(x,y): x \geq \beta   \text{ and } y\leq f(x,t)\bigg\}
\end{align*}
Below we list a number of specific grim reaper solutions that we use for convergence studies throughout the paper:
\small
\begin{itemize}
\item $(\theta_1,\theta_2,\theta_3)=(120^\circ,120^\circ,120^\circ)$\\
$(\sigma_{12},\sigma_{13},\sigma_{23})=(1,1,1)$\\
$\beta=\frac{1}{4}$
\[
  f(x,t) =
  \begin{cases}
    \frac{3}{2\pi}\log(\cos[\frac{2\pi}{3} x])-\frac{2\pi}{3} t & \text{if $0\leq x\leq\frac{1}{4}$} \\
    \frac{3}{2\pi}\log(\cos[\frac{2\pi}{3}(\frac{1}{2}-x)])-\frac{2\pi}{3} t  & \text{if $\frac{1}{4}<x\leq\frac{1}{2}$}
  \end{cases}
\]
\item $(\theta_1,\theta_2,\theta_3)=(90^\circ,135^\circ,135^\circ)$\\
$(\sigma_{12},\sigma_{13},\sigma_{23})=(1,1,\sqrt{2})$\\
$\beta=\frac{1}{4}$
\[
  f(x,t) =
  \begin{cases}
    \frac{1}{\pi}\log(\cos[\pi x])-\pi t & \text{if $0\leq x\leq\frac{1}{4}$} \\
    \frac{1}{\pi}\log(\cos[\pi(\frac{1}{2}-x)])-\pi t  & \text{if $\frac{1}{4}<x\leq\frac{1}{2}$}
  \end{cases}
\]

\item $(\theta_1,\theta_2,\theta_3)=(75^\circ,135^\circ,150^\circ)$\\ $(\sigma_{12},\sigma_{13},\sigma_{23})=\big(\frac{\sqrt{2}}{4}+\frac{\sqrt{6}}{4},\frac{\sqrt{2}}{2},\frac{1}{2}\big)$\\
$\beta=\frac{3}{46}(4\sqrt{2}-3)$
\[
\arraycolsep=0pt
  f(x,t) = \left\{
  \begin{array}{ll}
    \frac{1}{\frac{\pi}{6}(3+4\sqrt{2})}\log(\cos[ \frac{\pi}{6}(3+4\sqrt{2}) x])-\frac{\pi}{12}(3&+4\sqrt{2}) t\\
    &\text{if $0 \leq x\leq\frac{3}{46}(4\sqrt{2}-3)$} \\
    \frac{1}{\frac{\pi}{12}(3\sqrt{2}+8)}\log(\frac{2}{2^\frac{\sqrt{2}}{4}}\cos[\frac{\pi}{12}(3\sqrt{2}+8)(\frac{1}{2}-&x)])-\frac{\pi}{12}(3+4\sqrt{2}) t\\  &\text{if $\frac{3}{46}(4\sqrt{2}-3)<x\leq\frac{1}{2}$}
  \end{array} \right.
\]

\item $(\theta_1,\theta_2,\theta_3)=(75^\circ,135^\circ,150^\circ)$\\ $(\sigma_{12},\sigma_{13},\sigma_{23})=\big(\frac{\sqrt{2}}{4}+\frac{\sqrt{6}}{4},\frac{\sqrt{2}}{2},\frac{1}{2}\big)$\\
$\mu_{ij}=\frac{1}{\sigma_{ij}}$\\
$\beta=\frac{3}{14}$
\[
  f(x,t) =
  \begin{cases}
    \frac{6}{7\pi}\log(\cos[ \frac{7\pi}{6} x])-\frac{7\pi}{6} t & \text{if $0 \leq x\leq\frac{3}{14}$} \\
    \frac{6}{7\pi}\log(\sqrt{2}\cos[ \frac{7\pi}{6}(\frac{1}{2}-x)])-\frac{7\pi}{6} t  & \text{if $\frac{3}{14}<x\leq\frac{1}{2}$}
  \end{cases}
\]
\end{itemize}
\normalsize
\cref{fig:gr75} shows the exact solutions at time $t=0$ and $t=\frac{18}{512}$. In the case where $\theta_2=\theta_3$, $\Sigma_2(t)$ is the reflection of $\Sigma_3(t)$ around $x=\frac{1}{4}$, so we only need to track $\Sigma_1(t)$ and $\Sigma_2(t)$ on $\big[0,\frac{1}{4}\big]\times \big[-\frac{1}{2},1\big]$. 
\subsection{The VIIM via Parameterized Curves}
We begin the section by describing two essential numerical procedures.
The first finds the distance between a parameterized curve and a given point, and the second evolves a parameterized curve by curvature motion.
We then discuss our implementation of the VIIM on the ``Grim Reaper'' test case using parameterized curves.  
\subsubsection{Distance Estimation for a Parameterized Curve}

Finding the distance between a given point and a parameterized curve is important for constructing the $\epsilon$ level sets and for the Voronoi reconstruction step. To measure this distance with high accuracy we interpolate the parameterized curve with not-a-knot cubic splines~\cite{DeBoor}. Let $\{(x_i,y_i)\}_{i=1}^n$ be points on a curve that are given as the graph of a function of $x$ and denote its cubic spline approximation as $\bar{f}(x)$. The cubic spline is a piecewise third order polynomial that is twice differentiable over $[x_1,x_n]$ with coordinates $y_i=\bar{f}(x_i)$. The signed distance between a given point $(\tilde{x},\tilde{y})$ and the point which is closest to it on the aforementioned curve, is \[\text{sgn}(\tilde{y}-\bar{f}(\tilde{x}))\min_{x \in [x_1,x_N]} \sqrt{(x-\tilde{x})^2+(\bar{f}(x)-\tilde{y})^2}.\]
We find the minimum using Newton's method (\cref{alg:dist}).
\begin{algorithm}
  \caption{Netwon's method for finding distance between a parameterized curve and a point}
    \label{alg:dist}
  \begin{algorithmic}[1]
  \STATE Given a point $(\tilde{x},\tilde{y})$, a cubic spline curve $y = \bar{f}(x)$, points $\{(x_i,y_i)\}_{i=1}^n$ with $y_i = \bar{f}(x_i)$ and tolerance $\delta$:
  	\STATE Let $I=\text{arg}\min_i (x_i-\tilde{x})^2+(y_i-\tilde{y})^2$  \STATE Set $x=x_I$ to be the starting point in Newton's method
  		\WHILE{$|(x-\tilde{x})+(\bar{f}(x)-\tilde{y})\bar{f}'(x)|\geq\delta$}
        \STATE{
        \begin{align*}
        x \gets x-\frac{(x-\tilde{x})+(\bar{f}(x)-\tilde{y})\bar{f}'(x)}{1+(\bar{f}'(x))^2+(\bar{f}(x)-\tilde{y})\bar{f}''(x)}
        \end{align*}}
        \ENDWHILE
        \STATE Output $\text{sgn}(\tilde{y}-\bar{f}(\tilde{x}))\sqrt{(x-\tilde{x})^2+(\bar{f}(x)-\tilde{y})^2}$.
  \end{algorithmic}
\end{algorithm}

Newton's method uses the first two derivatives of the cubic spline. The interpolant is known to converge to the true curve in $C^2$~\cite{cubicerror}.\\  


While the foregoing discussion applies to a curve that is a function of $x$, the above techniques can be use for any simple curve by applying an appropriate change of variables. In general, distance functions are not differentiable everywhere, however we only use the distance function at differentiable points in the ``Grim Reaper'' cases we consider.

\subsubsection{Curvature Motion for a Parameterized Curve}
Curvature motion for a parameterized curve, $\gamma(s,t) \, : \, [0,L] \times \mathbb{R}^+ \to \mathbb{R}^2$, is described by:
\begin{align}
\label{eq:curve}
&\gamma_t=\partial_s\bigg(\frac{\gamma_s}{|\gamma_s|}\bigg)\frac{1}{|\gamma_s|}
\end{align}

The differential equation (\ref{eq:curve}) is implemented by a fully implicit Euler scheme where each iteration involves solving multiple tridiagonal linear systems.
The scheme in time and space is
\begin{align*}
&-\frac{ \delta t }{h^{k+1}_i h^{k+1}_{i-\frac{1}{2}}}\gamma^{k+1}_{i-1}+\bigg[1+\frac{ \delta t }{h^{k+1}_i}\bigg( \frac{1}{h^{k+1}_{i-\frac{1}{2}}}+ \frac{1}{h^{k+1}_{i+\frac{1}{2}}}\bigg)\bigg]\gamma^{k+1}_i  - \frac{ \delta t }{h^{k+1}_i h^{k+1}_{i+\frac{1}{2}}}\gamma^{k+1}_{i+1} =\gamma^{k}_i 
\end{align*}
where $h^k_i=\frac{1}{2}|\gamma^k_{i-1}-\gamma^k_{i+1}|$, $h^k_{i+\frac{1}{2}}=|\gamma^k_{i}-\gamma^k_{i+1}|$ and $\delta t$ is the time step size for the finite difference scheme (different than $\sigma\Delta t$, the total time we evolve the curve between reconstructions).
At each step we use the Newton iteration
\begin{align*}
&-\frac{ \delta t }{h^{k(l)}_i h^{k(l)}_{i-\frac{1}{2}}}\gamma^{k(l+1)}_{i-1}+\bigg[1+\frac{ \delta t }{h^{k(l)}_i}\bigg( \frac{1}{h^{k(l)}_{i-\frac{1}{2}}}+ \frac{1}{h^{k(l)}_{i+\frac{1}{2}}}\bigg)\bigg]\gamma^{k(l+1)}_i  - \frac{ \delta t }{h^{k(l)}_i h^{k(l)}_{i+\frac{1}{2}}}\gamma^{k(l+1)}_{i+1} =\gamma^{k}_i 
\end{align*}
until $|\sum_{i=0}^N h_i^{k(l)}-\sum_{i=0}^N h_i^{k(l+1)}| < \delta$ for a small $\delta>0$ and then set $\gamma^{k+1}=\gamma^{k(l)}$.\\

We detail how we handle all the boundary conditions in the $\theta_2=\theta_3$ case in \cref{tab:BC}. Then $\gamma(s,\sigma \Delta t)$ is given by $\gamma^K$ for $K=\frac{\sigma \Delta t}{\delta t}$ with initial value $\gamma^0=\gamma(s,0)$. Note that the choice of $\delta t$ is independent from $\Delta t$. In our numerical studies we choose $\delta t$ so small that the contribution to the overall error from the numerical solution $\gamma(s,\sigma \Delta t)$ is negligible.
\begin{table}
\begin{center}
\caption{Boundary Conditions for PDE \cref{eq:curve}}
\label{tab:BC}
\begin{tabular}{|c|c|c|}
\hline
Case&$x$ Boundary Condition&$y$ Boundary Condition\\
\hline
$x(0)=0$&$x(-s)=-x(s)$&$y(-s)=y(s)$\\
\hline
$x(L)=.25$&$x(L+s)=.5-x(L-s)$&$y(L+s)=y(L-s)$\\
\hline
$y(L)=-.5$&$x(L+s)=x(L-s)$&$y(L+s)=-1-y(L-s)$\\
\hline

\end{tabular}
\end{center}
\end{table}

\subsubsection{The Implementation of the VIIM using Parameterized Curves}

With these two tools in hand we can implement the VIIM using parameterized curves.
At every iterate we track a series of $(x,y)$ points that parameterize the interface of the sets. We will first consider the the case where $\theta_2=\theta_3$ (see the first two images in \cref{fig:gr75}). Since $\Sigma_3^k$ is the reflection of $\Sigma_2^k$ around $x=.25$ we only need to track the interface $\Gamma_{12}$. The interface $\Gamma_{12}$ remains a function of $x$ at every time step. Thus to parameterize $\Gamma_{12}$, we fix $x$ values in $[0,.25]$ and update corresponding $y$ values for each step in the VIIM. The simulation of mean curvature motion for time $T$ using parameterized curves is detailed below:
\begin{enumerate}
\item Choose $\Delta t$, the time between reconstructions, and $n$, the number of points. Set $N=T/\Delta t$.
\item Pick $\{x_i\}_{i=1}^n \in [0,.25]$ and set $y_i^0=f(x_i,0)$. In our implementation $x_i$'s are chosen so that $(x_{i-1}-x_i)^2+(f(x_{i-1},0)-f(x_{i},0))^2$ are all equal for $i=2,3,\ldots,n$.
\item For $k=1,\ldots,N$ do the following steps:
\begin{enumerate}
\item Build parameterization $\gamma^{\epsilon+}=\{(x,y):d_{\Sigma_1^{k-1}}(x,y)=\epsilon\}$:  For each $x_i$, we find $y^{\epsilon+}_i$ such that $d_{\Sigma_1^{k-1}}(x_i,y^{\epsilon+}_i)=\epsilon$, then $(x_i,y^{\epsilon+}_i)$ is a parameterization of $\gamma^{\epsilon+}$. To find  $d_{\Sigma_1^{k-1}}(x_i,y^{\epsilon+}_i)=\epsilon$, we use the secant method 
\begin{align}
\label{eq:secant}
\bar{y}^{n} \gets \bar{y}^{n-1}+(d_{\Sigma_1^{k-1}}(x_i,\bar{y}^{n-1})-\epsilon)\Bigg(\frac{\bar{y}^{n-1}-\bar{y}^{n-2}}{d_{\Sigma_1^{k-1}}(x_i,\bar{y}^{n-1})-d_{\Sigma_1^{k-1}}(x_i,\bar{y}^{n-2})}\Bigg)
\end{align}
until \begin{equation}
\label{eq:stop}
|d_{\Sigma_1^{k-1}}(x_i,\bar{y}^n)-\epsilon|<\delta
\end{equation}
for $\delta>0$. In our tests $|\frac{d}{dy} d_{\Sigma_1^{k-1}}(x_i,y)|>\frac{1}{4}$ near the true solution, so \cref{eq:stop} can be used to bound the error in $\bar{y}^n$. Similar statements hold when we employ the secant method in subsequent steps. We choose $\delta$ small enough so that the error from estimating $\bar{y}^n$ is negligible compared to the overall error. 
\item Build parameterization $\gamma^{\epsilon-}=\{(x,y):d_{\Sigma_2^{k-1}}(x,y)=\epsilon\}$: We let $x^{\epsilon-}_i=x_i\frac{.25-\epsilon}{.25}$ and find $y^{\epsilon-}_i$ such that $d_{\Sigma_2^{k-1}}(x^{\epsilon-}_i,y^{\epsilon-}_i)=\epsilon$ using the secant method. Then
\[
\gamma_i^{\epsilon-}=\begin{cases}
   (x^{\epsilon-}_i,y^{\epsilon-}_i) & \text{if $i\leq n$} \\
    \bigg(x^{\epsilon-}_n,y^{\epsilon-}_n+(i-n)\frac{(-.5-y^{\epsilon-}_N)}{M}\bigg)  & \text{if $n+1\leq i \leq n+m$}
  \end{cases}
\]
for some choice of $m$. We choose $m$ such that \[ \frac{(-.5-y^{\epsilon-}_n)}{m}\approx \sqrt{(x^{\epsilon-}_{n-1}-x^{\epsilon-}_{n})^2+(y^{\epsilon-}_{n-1}-y^{\epsilon-}_{n})^2}.\]
\item Evolve $\gamma^{\epsilon+}$ and $\gamma^{\epsilon-}$ by \cref{eq:curve} for time $\sigma_1\Delta t$ and $\sigma_2\Delta t$ respectively. Now $\gamma^{\epsilon+}=\partial \Sigma_1^{k-1}(\sigma_1 \Delta t,\epsilon)$ and $\gamma^{\epsilon-}=\partial \Sigma_2^{k-1}(\sigma_2 \Delta t,\epsilon)$.
\item  For each $x_i$ find $\tilde{y}_i$ such that
    \begin{align*}
    d_{\Sigma_1^{k-1}(\sigma_1 \Delta t,\epsilon)}(x_i,\tilde{y}_i)=d_{\Sigma_2^{k-1}(\sigma_2 \Delta t,\epsilon)}(x_i,\tilde{y}_i).
    \end{align*}
We find $\tilde{y}_i$ using the secant method (the update is similar to \cref{eq:secant}). To use $\gamma^{\epsilon-}$ to calculate $d_{\Sigma_2(\sigma_2 \Delta t,\epsilon)}(x,y)$, we apply a change of coordinates to make $\gamma^{\epsilon-}$ a function of $x$.
    \item Then assign each $\tilde{y}_i$ to $y^{k}_i$.  
\end{enumerate}
\item Then $\{(x_i,y^N_i)\}_{i=1}^n$ gives a parameterization of the interface $\Gamma_{12}$ at time $T$.
\end{enumerate}

\subsection{Experimental Results}
In experiments with the VIIM, it suffices to focus on the symmetric cases $\theta_2 = \theta_3$ to demonstrate non-convergence.
Denote the output of the VIIM at time $T$ by $\widehat{\Sigma}_j(T)$ and the true solution by $\Sigma_j(T)$.
We track essentially the relative error (RE) of the area of symmetric difference in phase $\Sigma_1$:
 \begin{align}
 \label{rel}
\frac{\left| \widehat{\Sigma}_1(T) \, \triangle \, \Sigma_1(T) \right|}{\Big| \Sigma_1(T) \, \triangle \, \Sigma_1(0) \Big|}
 \end{align}
but restricted to $\{ (x,y) \, : 0 \leq x \leq 0.21 \}$ to exclude a small neighborhood around the junction at $x=\frac{1}{4}$.

Each of the simulations uses the following parameters
\begin{itemize}
\item The total time the system is evolved: $T=18/512$.
\item Number of points tracked on the parameterized curve: $n=2048$.
\item Step size in \cref{eq:curve}: $\delta t=2^{-12}\sigma\Delta t$.
\item In the equal surface tension case ($\theta_1=120^\circ$): $\sigma_1=\sigma_2=1$
\item We use $\theta_1=90^\circ$ in our test for the unequal surface tension case, so that $\sigma_1=2-\sqrt{2}$ and $\sigma_2=\sqrt{2}$
\end{itemize}
Refining in $\delta t$ or $n$ did not significantly change the relative error. Additionally, we simulate the formal limit, $\displaystyle \lim_{\epsilon \to 0^+} \displaystyle  \lim_{\Delta t \to 0}$, by setting $\epsilon\propto\Delta t^{1/4}$.\\

The results are collected in \cref{tab:120a} through \cref{tab:90c}.
 In none of the limit cases does the unequal surface tension case converge to the correct curve. This is seen in the non-vanishing relative error for $\theta_1 \neq 120^\circ$. Later we will give an alternative to the Voronoi reconstruction step that, in our numerical tests, convergences in the unequal surface case.
\begin{table}
\parbox{.49\linewidth}{
\begin{center}
\caption{$\lim_{\epsilon \to 0^+}\lim_{\Delta t \to 0}$\\  $\theta_1=120^\circ$}
\label{tab:120a}
\begin{tabular}{|c|c|c|c|}
\hline
$\Delta t$&$\epsilon$&RE &Order\\
\hline
$2^{-13}$&$2^{-28/4}$&0.0339&-\\
\hline
$2^{-14}$&$2^{-29/4}$&0.0269&0.33\\
\hline
$2^{-15}$&$2^{-30/4}$&0.0214&0.33\\
\hline
$2^{-16}$&$2^{-31/4}$&0.0172&0.32\\
\hline
\end{tabular}
\end{center}
}
\parbox{.49\linewidth}{
\begin{center}
\caption{$\lim_{\epsilon \to 0^+}\lim_{\Delta t \to 0}$\\ $\theta_1=90^\circ$   }
\label{tab:90a}
\begin{tabular}{|c|c|c|c|}
\hline
$\Delta t$&$\epsilon$&RE &Order\\
\hline
$2^{-13}$&$2^{-28/4}$&0.0361&-\\
\hline
$2^{-14}$&$2^{-29/4}$&0.0517&-\\
\hline
$2^{-15}$&$2^{-30/4}$&0.0667&-\\
\hline
$2^{-16}$&$2^{-31/4}$&0.0813&-\\
\hline
\end{tabular}
\end{center}
}
\end{table}

\begin{table}
\caption{$\lim_{\epsilon=c \sqrt{\Delta t} \to 0^+}$\\$\theta_1=120^\circ$}
\label{tab:120b}
\begin{center}
\begin{tabular}{|c|c|c|c|c|}
\hline
$\Delta t$&RE: $c=4$&Order&RE: $c=2$&Order\\
\hline
$2^{-13}$&0.1675&-&0.0833&-\\
\hline
$2^{-14}$&0.1078&0.636&0.0561&0.572\\
\hline
$2^{-15}$&0.0714&0.595&0.0383&0.551\\
\hline
$2^{-16}$&0.0482&0.566&0.0264&0.537\\
\hline
\end{tabular}
\end{center}
\end{table}

\begin{table}
\caption{$\lim_{\epsilon=c \sqrt{\Delta t} \to 0^+}$ \\ $\theta_1=90^\circ$}
\label{tab:90b}
\begin{center}
\begin{tabular}{|c|c|c|c|c|}
\hline
$\Delta t$&RE: $c=4$&Order&RE: $c=2$&Order\\
\hline
$2^{-13}$&0.0256&-&0.0524&-\\
\hline
$2^{-14}$&0.0826&-&0.0791&-\\
\hline
$2^{-15}$&0.1166&-&0.0963&-\\
\hline
$2^{-16}$&0.1379&-&0.1077&-\\
\hline
\end{tabular}
\end{center}
\end{table}

\begin{table}
\parbox{.49\linewidth}{
\begin{center}
\caption{$\lim_{\Delta t \to 0}$ with $\epsilon=0$\\ $\theta_1=120^\circ$}

\label{tab:120c}
\begin{tabular}{|c|c|c|}
\hline
$\Delta t$&RE &Order\\
\hline
$2^{-13}$&0.0071&-\\
\hline
$2^{-14}$&0.0050&0.51\\
\hline
$2^{-15}$&0.0035&0.51\\
\hline
$2^{-16}$&0.0025&0.50\\
\hline
$2^{-17}$&0.0017&0.50\\
\hline
$2^{-18}$&0.0012&0.50\\
\hline
$2^{-19}$&0.0009&0.50\\
\hline
\end{tabular}
\end{center}
}
\parbox{.49\linewidth}{
\begin{center}
\caption{$\lim_{\Delta t \to 0}$ with $\epsilon=0$\\ $\theta_1=90^\circ$}
\label{tab:90c}
\begin{tabular}{|c|c|c|}
\hline
$\Delta t$&RE &Order\\
\hline
$2^{-13}$&0.0056&-\\
\hline
$2^{-14}$&0.0065&-\\
\hline
$2^{-15}$&0.0071&-\\
\hline
$2^{-16}$&0.0075&-\\
\hline
$2^{-17}$&0.0077&-\\
\hline
$2^{-18}$&0.0078&-\\
\hline
$2^{-19}$&0.0079&-\\
\hline
\end{tabular}
\end{center}
}
\end{table}

\section{Threshold Dynamics}
\label{sec:thres}
In this section we present convergence studies for the threshold dynamics algorithm of \cite{esedog2015threshold} using the parametrized curve implementation developed above.
There is by now ample evidence, including a conditional proof ~\cite{lauxotto}, for the convergence of this algorithm to the correct limit, including very general unequal surface tensions cases.
This section is thus meant as a verification of the parametrized curve implementation (rather than threshold dynamics, which is not in doubt), and give confidence to the non-convergence results it yielded on the VIIM, presented in the previous section.

\begin{table}
\parbox{.49\linewidth}{

\begin{center}
\caption{Threshold Dynamics\\ $\theta_1=120^\circ$}
\label{tab:TD120}
\begin{tabular}{|c|c|c|c|}
\hline
$\Delta t$&$n$&RE &Order\\
\hline
$2^{-13}$&512&0.0081&-\\
\hline
$2^{-14}$&1024&0.0056&0.523\\
\hline
$2^{-15}$&2048&0.0039&0.515\\
\hline
$2^{-16}$&4096&0.0028&0.510\\
\hline
$2^{-17}$&8192&0.0019&0.507\\
\hline
$2^{-18}$&16384&0.0014&0.504\\
\hline
$2^{-19}$&32768&0.0010&0.504\\
\hline
\end{tabular}
\end{center}
}
\parbox{.49\linewidth}{

\begin{center}
\caption{Threshold Dynamics\\
$\theta_1=90^\circ$}
\label{tab:TD90}
\begin{tabular}{|c|c|c|c|}
\hline
$\Delta t$&$n$&RE&Order\\
\hline
$2^{-13}$&512&0.0095&-\\
\hline
$2^{-14}$&1024&0.0067&0.516\\
\hline
$2^{-15}$&2048&0.0047&0.511\\
\hline
$2^{-16}$&4096&0.0033&0.507\\
\hline
$2^{-17}$&8192&0.0023&0.505\\
\hline
$2^{-18}$&16384&0.0016&0.505\\
\hline
$2^{-19}$&32768&0.0012&0.501\\
\hline
\end{tabular}
\end{center}
}
\end{table}
 We use the following parameters:
 \begin{itemize}
 \item Each of the following simulations is evolved for time $T=18/512$.
 \item For the equal surface tension case $\theta_1=120^\circ$, we use $\sigma_{12}=\sigma_{13}=\sigma_{23}=1$.
 \item For the unequal surface tension case with $\theta_1=90^\circ$, we use $\sigma_{12}=\sigma_{13}=1$ and $\sigma_{23}=\sqrt{2}$.
 \end{itemize}

The results are in \cref{tab:TD120} and \cref{tab:TD90}. 
We see convergence to the correct solution, including in the unequal surface tension case, bolstering our confidence in the algorithm and the parametrized curve implementation developed and used in this paper.
It is thus highly unlikely that the non-convergence observed in the previous section with the VIIM is due to the parametrized curve representation; it is likely due to the VIIM itself.

\section{Correcting the VIIM: Dictionary Mapping}
\label{sec:alg}
Before the variational formulation of threshold dynamics given in \cite{esedog2015threshold} extended the original algorithm of \cite{merriman1994motion} from equal to arbitrary surface tensions in a systematic manner, a more heuristic extension was proposed by Ruuth in \cite{ruuth}.
In this approach, a projection step is employed to ``force'' the correct Herring angle conditions at any triple junction, while multiple ($\geq 4$) junctions are treated more heuristically.
Ruuth's projection is designed so that the stationary configuration for the underlying curvature flow of three flat interfaces meeting at a triple junction with the correct Herring angles remains {\it fixed} under one iteration of the overall algorithm.
Motivated by Ruuth's approach, in this section we propose a new algorithm: the dictionary mapping implicit interface method (DMIIM), which replaces the Voronoi reconstruction step of the VIIM with a {\it dictionary reconstruction step} that is designed to have as a fixed point three flat interfaces meeting with the correct Herring angles. The three phases in such a configuration consist of sectors and after evolution by curvature motion we want to restore these phases to their original form by our reconstruction. As such our dictionary reconstruction step is based on the curvature evolution of sectors.
A heuristic extension to arbitrary number of phases, much as in \cite{ruuth}, is also discussed.
While an analogue of the more systematic approach of \cite{esedog2015threshold} would be more satisfactory, no such variational formulation for the VIIM (that is simple and efficient to implement) is currently available -- a matter that remains under investigation.\\

We describe dictionary reconstruction in the case of three phases.
The first step is to build a {\it template map} for a triple junction of interfaces with surface tensions $\sigma_1$, $\sigma_2$, and $\sigma_3$.
To that end, let $(\theta_1,\theta_2,\theta_3)$ be the triple junction angles corresponding to these surface tensions as given by equation \cref{herring}, and fix a time step size $\Delta t$.
Let $\Omega_i(0)$ for $i\in\{1,2,3\}$ denote the sector
\[ \Omega_i(0) = \Big\{ (\theta,r) \, : r\geq 0 \mbox{ and } \sum_{j = 0}^{i-1} \theta_j \leq \theta \leq \sum_{j=0}^{i} \theta_j \Big\} \]
in polar coordinates, with the proviso $\theta_0=0$.
Let $\Omega_i(\sigma_i\Delta t)$ be the evolution of $\Omega_i(0)$ via motion by mean curvature for time $\sigma_i \Delta t$. Recall that $d_{\Omega_i(\sigma_i\Delta t)}$ is the signed distance function to $\Omega_i(\sigma_i\Delta t)$, see \cref{eq:distance}. We will define the {\it template map}, $\Phi:\mathbb{R}^2 \to \mathbb{R}^3$, as 

\begin{align}
\label{eq:ts}
\Phi(x,y)=(d_{\Omega_1(\sigma_1\Delta t)}(x,y),d_{\Omega_2(\sigma_2\Delta t)}(x,y),d_{\Omega_3(\sigma_3\Delta t)}(x,y))
\end{align}
and define the {\it template surface} $S \subset \mathbb{R}^3$ as the image of $\mathbb{R}^2$ under $\Phi$. The template map, $\Phi$, maps points in $\mathbb{R}^2$ to distances to the evolved sectors. The map $\Phi$ is injective and we will use $\Phi^{-1}:S \to \mathbb{R}^2$ in our reconstruction algorithm. \\

\begin{figure}

  \begin{center}
  \subfloat{\includegraphics[width=1\textwidth]{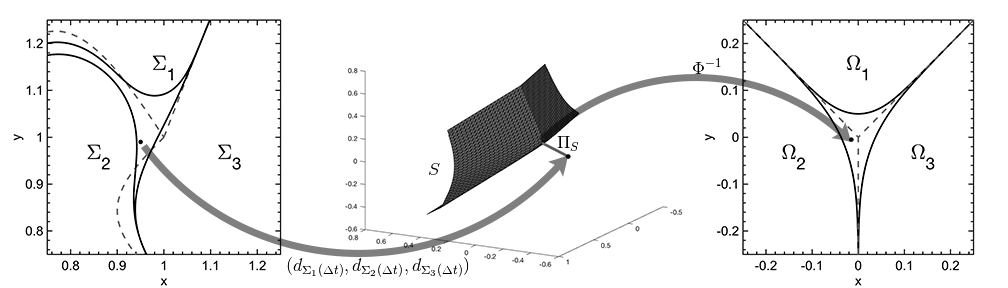}}

\caption{\footnotesize A schematic of the dictionary reconstruction step. The dashed lines are the interfaces at $T=0$ and the solid lines are sets at time $T=\Delta t$. In this example the point would be allocated to $\Sigma_2$.}
\label{fig:S}
  \end{center}
\end{figure}

Recall in the VIIM that $\Sigma_i^k$ is phase $i$ at time step $k$. We describe how the DMIIM reconstructs the new phases, $\Sigma_1^{k+1}$, $\Sigma_2^{k+1}$, and $\Sigma_3^{k+1}$ from $\Sigma_1^{k}(\sigma_1 \Delta t)$, $\Sigma_2^{k}(\sigma_2 \Delta t)$, and $\Sigma_3^{k}(\sigma_3 \Delta t)$. Define $\Pi_{S}: \mathbb{R}^3 \to \mathbb{R}^3$ as the closest point projection onto $S$ (with respect to the standard Euclidean distance in $\mathbb{R}^3$).
We define the reconstructed phases at the $(k+1)$-st time step as
\begin{align}
\label{eq:proj} 
\Sigma^{k+1}_i= \Big\{ \mathbf{z}\in D \, : \,  \Phi^{-1} \circ \Pi_{S} \big( d_{\Sigma^k_1(\sigma_1 \Delta t)}(\mathbf{z}),d_{\Sigma^k_2(\sigma_2\Delta t)}(\mathbf{z}),d_{\Sigma^k_3(\sigma_3\Delta t)}(\mathbf{z}) \big) \in \Omega_i(0) \Big\}.
\end{align}
We are thus assigning the point $\mathbf{z}$ to the phase (or phases) whose corresponding sector contains the preimage of
\[\Pi_{S}(d_{\Sigma_1(\sigma_1\Delta t)}(\mathbf{z}),d_{\Sigma_2(\sigma_2\Delta t)}(\mathbf{z}),d_{\Sigma_3(\sigma_3\Delta t)}(\mathbf{z})) \in S\]
under the one-to-one map $\Phi$.
The configuration $\big( \Omega_1(0) , \Omega_2(0) , \Omega_3(0) \big)$ is clearly fixed under the DMIIM algorithm, which thus treats a triple junction with the correct Herring angles and straight interfaces exactly as it should.
(Note that having these triple junctions fixed is a reassuring but not necessary condition for convergence of such algorithms; e.g. threshold dynamics \cite{esedog2015threshold} does not have this property in general). See \cref{fig:S} for a schematic of dictionary reconstruction.\\

Given the surface tensions $(\sigma_1,\sigma_2,\sigma_3)$ and a time step size $\Delta t > 0$, the corresponding projection surface $S$ is neither bounded nor smooth.
Moreover, at its non-empty ridge, the closest point projection map $\Pi_S \, : \, \mathbb{R}^3 \to S$ contains multiple points.
However, $\Pi_S(\mathbf{x})$ is non-empty at any $\mathbf{x}\in\mathbb{R}^3$.\\

The following two claims ensure that $\Pi_S$ is well defined (i.e. $\Pi_S(\mathbf{x})$ contains a single point) and smooth in a neighborhood of the vacuum and overlap regions formed by evolving each phase by mean curvature motion, starting from an initial configuration of three smooth curves meeting at almost the correct Herring angles. The claims are straightforward but tedious to check, so the proofs are omitted.\\

The first concerns the surface $S$ and follows from properties of the self-similar solution of curvature motion discussed in \cref{SSSCM}:
\begin{claim}
Given $\mathbf{r}=(r_1,r_2,r_3) \in (\mathbb{R}^+)^3$, let $T_{\mathbf{r}} = \cap_{j=1}^3 \{ \mathbf{x} \in \mathbb{R}^2 \, : \, d_{\Omega_j(0)}(\mathbf{x}) < r_j\}$ denote a neighborhood of the stationary triple junction with the exact Herring angles.
There exists $\mathbf{r}>\mathbf{0}$ and $\varepsilon > 0$ such that the closest point projection map $\Pi_S \, : \, \mathbb{R}^3 \to S$ is well defined and smooth on $N_\varepsilon = \{ \mathbf{x} \in \mathbb{R}^3 \, : \, d(\mathbf{x},\Phi(T_{\mathbf{r}})) < \varepsilon \}$.
\end{claim}
The second can be checked e.g. using the comparison principle satisfied by motion by mean curvature:
\begin{claim}
\label{claim:cond}
Let $\mathbf{r}, \varepsilon, T_{\mathbf{r}}$, and $N_\varepsilon$ be as in Claim 1.
Let $(\partial\Sigma_i(0)) \cap (\partial\Sigma_j(0))$, $i \not= j$, be smooth curves meeting at a triple junction with angles $(\theta'_1,\theta'_2,\theta'_3)$.
There exist $\delta > 0$ and $T>0$ such that if $||(\theta_1,\theta_2,\theta_3) - (\theta'_1,\theta'_2,\theta'_3)|| \leq \delta$ and $\Delta t \leq T$, then
$$ \frac{1}{\sqrt{\Delta t}}( d_{\Sigma_1(\sigma_1 \Delta t)}(\mathbf{x}) , d_{\Sigma_2(\sigma_2 \Delta t)}(\mathbf{x}) , d_{\Sigma_3(\sigma_3 \Delta t)}(\mathbf{x}) ) \in N_{\varepsilon} \mbox{ for any } \mathbf{x} \in \bigcup_{\substack{ 0 \leq t\leq \Delta t\\ j \in \{1,2,3\}}} \partial \Sigma_j(\sigma_j t). $$
\end{claim}

We next explain one way to extend the DMIIM algorithm just described for three phases to the $n$-phase setting.
Let $\sigma_1,\sigma_2,\ldots,\sigma_n$ be the surface tensions associated with the $n$ phases.
Let $\mathbf{I}(\mathbf{z},\cdot) \, : \, \{ 1 , 2 , \ldots , n \} \to \{ 1 , 2 , \ldots , n \}$ be a bijection (permutation of the indices) so that
\begin{equation*}
d_{\Sigma^k_{\mathbf{I}(\mathbf{z},1)}(\sigma_{\mathbf{I}(\mathbf{z},1)}\Delta t)}(\mathbf{z}) \geq d_{\Sigma^k_{\mathbf{I}(\mathbf{z},2)}(\sigma_{\mathbf{I}(\mathbf{z},2)}\Delta t)}(\mathbf{z}) \geq \cdots \geq
d_{\Sigma^k_{\mathbf{I}(\mathbf{z},n)}(\sigma_{\mathbf{I}(\mathbf{z},n)}\Delta t)}(\mathbf{z})
\end{equation*}
so that $\mathbf{I}(\mathbf{z},j)$ is the label of the phase with the $j$-th largest signed distance function at the point $\mathbf{z}$.
Our simple extension, similar to the one in \cite{ruuth}, is to allocate each $\mathbf{z}$ using the very same dictionary mapping discussed above, where the three phases used in the construction of the projection surface are the closest three to $\mathbf{z}$.
Under this rule the new sets become
\small
\begin{multline}
\label{eq:ext}
\Sigma^{k+1}_j = \Big\{ \mathbf{z} \in D \, : \, \mathbf{I}^{-1}(\mathbf{z},j) \leq 3 \mbox{ and }\\
\Phi^{-1} \circ \Pi_S \big( d_{\Sigma^k_{\mathbf{I}(\mathbf{z},1)}(\sigma_{\mathbf{I}(\mathbf{z},1)}\Delta t)}(\mathbf{z}),d_{\Sigma^k_{\mathbf{I}(\mathbf{z},2)}(\sigma_{\mathbf{I}(\mathbf{z},2)}\Delta t)}(\mathbf{z}),d_{\Sigma^k_{\mathbf{I}(\mathbf{z},3)}(\sigma_{\mathbf{I}(\mathbf{z},3)}\Delta t)}(\mathbf{z}) \big) \in \Omega _{\mathbf{I}^{-1}(\mathbf{z},j)}(0) \Big\}
\end{multline}
\normalsize
where $\Omega_j$ for $j\in\{1,2,3\}$, the map $\Phi$ and the surface $S$ are constructed at each $\mathbf{z}\in D$ as in the three phase case, using the surface tensions $\sigma_{\mathbf{I}(\mathbf{z},1)}$, $\sigma_{\mathbf{I}(\mathbf{z},2)}$, and $\sigma_{\mathbf{I}(\mathbf{z},3)}$. High degree junctions are thus treated heuristically by this method. Indeed, it is not hard to come up with alternatives to the simple extension \cref{eq:ext} explained above. Although we expect all such natural extensions to behave mostly the same, subtle differences between them cannot be ruled out at this point. We took the simplest example \cref{eq:ext}, and while points near high degree junctions can be far away from the template surface, we will show that it behaves reasonably in \cref{sec:four}.\\

The DMIIM also allows some control over mobilities. Each phase $\Sigma_i$ can be assigned a mobility $\mu_i$. In the evolution step the level sets are evolved for time $\mu_i \sigma_i \Delta t$ by \ref{eq:ls} and we use $\Omega_i(\mu_i \sigma_i \Delta t)$ in the construction of the template surface. The mobility at the interface $\Gamma_{ij}$ becomes
\begin{align}
\label{eq:mobile}
\mu_{ij}=\frac{\mu_i \sigma_i+\mu_j \sigma_j}{2\sigma_{ij}}.
\end{align}
With this extra flexibility, the DMIIM allows the specification of the $\binom{n}{2}$ physically relevant surface tensions at the interfaces, $\Gamma_{ij}$, and the products, $\mu_i\sigma_i$, for each phase. The mobility at the interface $\Gamma_{ij}$ is constrained by \cref{eq:mobile}. We give an example in \cref{sec:four}.\\

\cref{alg:DMIIM} details the steps in the DMIIM.  In the next two subsections we describe how to find $\Omega_i(t)$ with high accuracy and how the projection is preformed. \\


\begin{algorithm}
  \caption{Dictionary Mapping Implicit Interface Method
    \label{alg:DMIIM}}
  \begin{algorithmic}[1]
  	\STATE Given $\Sigma^0_1,\Sigma^0_2,\ldots,\Sigma^0_n$, $\sigma_1,\sigma_2,\ldots,\sigma_n$, $\mu_1,\mu_2,\ldots,\mu_n$, $\Delta t$, and $T$.
    \STATE Let $N=T/\Delta t$. Define the {\em reduced mobilities} $\bar{\mu}_i = \mu_i \sigma_i$.
      \FOR{$k \gets 1 \textrm{ to } N$}
          \STATE Evolve each $\Sigma_i^{k}(0)$ by time $\bar{\mu}_i \Delta t$ to get $\Sigma_i^{k}(\bar{\mu}_i\Delta t)$.
          \STATE Construct the new phases
          \begin{multline*}
\Sigma^{k+1}_i= \Big\{ \mathbf{z} \in D \, : \, \mathbf{I}^{-1}(\mathbf{z},i) \leq 3 \mbox{ and }\\
\Phi^{-1}\circ \Pi_S \big( d_{\Sigma^k_{\mathbf{I}(\mathbf{z},1)}(\bar{\mu}_{\mathbf{I}(\mathbf{z},1)} \Delta t)}(\mathbf{z}),d_{\Sigma^k_{\mathbf{I}(\mathbf{z},2)}(\bar{\mu}_{\mathbf{I}(\mathbf{z},2)} \Delta t)}(\mathbf{z}),d_{\Sigma^k_{\mathbf{I}(\mathbf{z},3)}(\bar{\mu}_{\mathbf{I}(\mathbf{z},3)} \Delta t)}(\mathbf{z}) \big)\\ \in \Omega _{\mathbf{I}^{-1}(\mathbf{z},i)}(0) \Big\}
\end{multline*}
      \ENDFOR
  \end{algorithmic}
\end{algorithm}

\subsection{Self-Similar Solution to Curvature Motion}
\label{SSSCM}
Building the template map in the DMIIM algorithm requires precomputing the solution of motion by mean curvature of a sector (denoted $\Omega(t)$ in the previous section) to high accuracy. The set $\Omega(t)$ is related to a self-similar solution of mean curvature motion; thus the computation of $\Omega(t)$ can be reduced to the following ODE:
\begin{align}
\label{eq:ode}
\begin{split}
&\phi''(x)=\frac{1}{2}(\phi(x)-x\phi'(x))(1+(\phi'(x))^2)\\
&\phi'(0)=0\\
&\lim_{x \to \infty}\phi(x)=\infty\\
&\phi(0)=\phi_0>0\\
\end{split}
\end{align}
on the domain $x\in (-\infty,\infty)$ and $\phi$ is even. Instead of the last condition of \cref{eq:ode}, we could specify a $M$ such that $\displaystyle \lim_{x \to \infty}\phi'(x)=M>0$. There is a bijective map between $\phi_0$ and $M$~\cite{ishimura1995curvature}.
We next explain how \cref{eq:ode} arises.\\

For a curve given as the graph of a function $u(x,t)$, motion by curvature is described by the PDE
\[
u_t=\frac{u_{xx}}{1+u_{x}^2}.
\]
When $u(x,t=0)=M|x|$ for some $M$ then \[u(x,t)=
\sqrt{t}\phi(x/\sqrt{t})\] where $\phi$ satisfies (\ref{eq:ode}). Let the positive y-axis bisect the sector $\Omega(0)$ with angle $\theta$, then $\partial\Omega(0)=M|x|$ for $M=\cot(\theta/2)$ and \[\Omega(t)=\{(x,y):y\geq \sqrt{t}\phi(x/\sqrt{t})\}.\]
To find the numerical solution to the ODE (\ref{eq:ode}), we use the Newton iteration
\begin{align}
&\bigg[\frac{2}{h^2}-\frac{x_i}{2h}-\frac{x_i}{2h}\bigg(\frac{\phi_{i+1}^k-\phi_{i-1}^k}{2h}\bigg)^2\bigg]\phi_{i-1}^{k+1}-\bigg[\frac{4}{h^2}+1+\bigg(\frac{\phi_{i+1}^k-\phi_{i-1}^k}{2h}\bigg)^2\bigg]\phi_i^{k+1}\\&+\bigg[\frac{2}{h^2}+\frac{x_i}{2h}+\frac{x_i}{2h}\bigg(\frac{\phi_{i+1}^k-\phi_{i-1}^k}{2h}\bigg)^2\bigg]\phi_{i+1}^{k+1}=0
\end{align}
for $x_0=0$ and $x_N$ being sufficiently large. The boundary conditions are $\phi_{-1}=\phi_{1}$ and $\phi_N=Mx_N$. Below we prove that $\phi_N$ is close to $Mx_N$ as long as we choose $x_N$ large enough, justifying the second boundary condition. We offer an improved bound over 
\begin{align}
\label{eq:bound1}
\phi(x)=Mx+o\bigg(\frac{1}{x}\bigg),\, \text{as }x \rightarrow \infty.
\end{align}
given in~\cite{ishimura1995curvature}, where the ODE was previously studied.  The improved bound implies we do not need to take $x_N$ so large. In our simulations we choose $x_N$ to be 10. 

\begin{claim}
\label{thm:selfsol}
For a function $\phi(x)$ satisfying (\ref{eq:ode}) the following bounds hold:
\begin{align*}
|\phi(x)-xM|&\leq C_0e^{-\frac{x^2(1+M^2)}{4}},\,x\geq0\\
|\phi'(x)-M|&\leq \frac{C_0}{x}e^{-\frac{x^2(1+M^2)}{4}},\,x>0\\
|\phi''(x)|&\leq C_1e^{-\frac{x^2(1+M^2)}{4}},\,x\geq0
\end{align*}
where $C_0$ and $C_1$ only depend on $\phi_0$ (or $M$).
\end{claim}

We first need the following lemmas
\begin{lemma}
\label{lemma:ineq}
The function $\phi$ satisfies the following properties:

\begin{enumerate}
\item $\phi''(x)>0$.
\item $0\leq \phi'(x)<M$.
\item $\phi(x)>xM$.
\end{enumerate}
\end{lemma}
\begin{proof}
The proof of property 1 is given in~\cite{ishimura1995curvature}. As a result of $\phi''(x)>0$, $\phi'(x)$ is a strictly increasing function with $\displaystyle\lim_{x \to \infty}\phi'(x)=M$, so property 2 follows. To prove property 3, let $h(x)=\phi(x)-xM$ on $x\geq0$. By \cref{eq:bound1} $\displaystyle \lim_{x \to \infty} h(x)=0$ and property 2 implies
 $h'(x)=\phi'(x)-M<0$. Thus $h(x)>0$ for all $x\geq0$.
\end{proof}

\begin{lemma}
The function $\phi$ satisfies the first order differential equation
\[
\exp{\bigg(\frac{\phi^2(0)}{2}\bigg)}\phi^2(0)=(\phi(x)-x\phi'(x))^2\exp{\bigg(\frac{x^2}{2}+\frac{\phi^2(x)}{2}\bigg)}(1+(\phi'(x))^2)^{-1}.
\]
\end{lemma}
\begin{proof}
Rearrange 
\begin{align}
\label{eqn:cm}
\phi''(\tilde{x})=\frac{1}{2}(\phi(\tilde{x})-\tilde{x}\phi'(\tilde{x}))(1+(\phi'(\tilde{x}))^2)
\end{align}
to
\[
\frac{-2\tilde{x}\phi''(\tilde{x})}{\phi(\tilde{x})-\tilde{x}\phi'(\tilde{x})}=-\tilde{x}+\tilde{x}(\phi'(\tilde{x}))^2.
\]
By integrating both sides from 0 to $x$ we obtain
\begin{align}
\label{part1}
2\log{(\phi(x)-x\phi'(x))}-2\log{(\phi_0}=-\frac{1}{2}x^2-\int_0^x \tilde{x}(\phi'(\tilde{x}))^2 d\tilde{x}.
\end{align}
Now rearrange \cref{eqn:cm} to
\[
\frac{2\phi'(\tilde{x})\phi''(\tilde{x})}{1+(\phi'(\tilde{x}))^2}=\phi(\tilde{x})\phi'(\tilde{x})-\tilde{x}(\phi'(\tilde{x}))^2.
\]
By integrating both sides from 0 to $x$ we obtain

\begin{align}
\label{part2}
\log{(1+(\phi'(x))^2)}= \frac{1}{2}\phi^2(x)-\frac{1}{2}\phi^2_0 -\int_0^x \tilde{x}(\phi'(\tilde{x}))^2 d\tilde{x}.
\end{align}

Both \cref{part1} and \cref{part2} have a $-\int_0^x \tilde{x}\phi'(\tilde{x})^2 d\tilde{x}$ term. Solving for that term in \cref{part1} and \cref{part2} and setting the equations equal to each other results in
\[
\log{(1+(\phi'(x))^2)}-\frac{1}{2}\phi^2(x)+\frac{1}{2}\phi^2_0=2\log{(\phi(x)-x\phi'(x))}-2\log{(\phi_0)}+\frac{1}{2}x^2.
\]
The conclusion of the lemma follows from taking the exponential of both sides.
\end{proof}
Using the above two lemmas we can establish \cref{thm:selfsol}:
\begin{proof}(Of \cref{thm:selfsol})
Applying the lemmas we have
\begin{align*}
& (\phi(x)-x\phi'(x))^2\\
=&\exp{\bigg(-\frac{x^2}{2}+\frac{\phi^2_0}{2}-\frac{\phi^2(x)}{2}\bigg)}(1+(\phi'(x))^2)\phi^2_0\\
\leq&\exp{\bigg(-\frac{x^2}{2}-\frac{(Mx)^2}{2}\bigg)}(1+M^2)\phi^2_0\exp{\frac{\phi^2_0}{2}}
\end{align*}
Additionally, invoking inequalities 2 and 3 from \cref{lemma:ineq}.
\[|\phi(x)-x\phi'(x)|=|\phi(x)-xM+xM-x\phi'(x)|=|\phi(x)-xM|+|xM-x\phi'(x)|.\] The first inequality follows from observing that $|\phi(x)-xM|<|\phi(x)-x\phi'(x)|$ and the second from $|M-\phi'(x)|<\frac{1}{x}|\phi(x)-x\phi'(x)|$ for $x > 0$. Then using \cref{eqn:cm}  $|\phi''(x)|\leq\frac{1}{2}|\phi(x)-x\phi'(x)|(1+M^2)$ and the third inequality follows.
\end{proof}

\subsection{Projecting onto the Template Surface}
Another important step in the DMIIM algorithm introduced in \cref{sec:alg} is the closest point projection onto the template surface. Here, we discuss the details of a highly accurate method for projecting an arbitrary point $\mathbf{w}\in \mathbb{R}^3$ onto 
the template surface $S=\{\Phi(x,y):(x,y) \in \mathbb{R}^2\}$ for the template map $\Phi$ defined in \cref{eq:ts}. Define the function
\begin{align}
\label{eq:dist3d}
F(x,y)=\frac{1}{2}||\mathbf{w}-\Phi(x,y)||_2^2.
\end{align}
Denote $(x^*,y^*)$ as the minimum of $F$ or where $\Phi(x^*,y^*)=\Pi_S(\mathbf{w})$. 
To find $(x^*,y^*)$ we will use Newton's method, 
\[
\begin{bmatrix}
    x^{n+1}\\
    y^{n+1}
\end{bmatrix}=
\begin{bmatrix}
    x^{n}\\
    y^{n}
\end{bmatrix}-(D^2F)^{-1}\nabla F(x^n,y^n).
\]
To choose $(x^0,y^0)$, we make a coarse point cloud of the surface $S$ and choose $(x^0,y^0)$ so that $\Phi(x^0,y^0)$ is the nearest point to $\mathbf{w}$ in the point cloud.\\

The rest of this section details how to compute $\nabla F$ and $D^2F$. First we find partial derivatives of $F$ in terms of the distance functions, for example:
\begin{align*}
&\frac{\partial F}{\partial x}=\sum_{i=1}^3 \bigg(\frac{\partial}{\partial x}  d_{\Omega_i(\sigma_i\mu_i\Delta t)} \bigg)(d_{\Omega_i(\sigma_i\mu_i\Delta t)}(x,y)-\mathbf{w}_i)\\
&\frac{\partial^2 F}{\partial x^2}= \sum_{i=1}^3 \bigg(\frac{\partial^2 }{\partial x^2} d_{\Omega_i(\sigma_i\mu_i\Delta t)}\bigg)(d_{\Omega_i(\sigma_i\mu_i\Delta t)}-\mathbf{w}_i)+\bigg(\frac{\partial}{\partial x} d_{\Omega_i(\sigma_i\mu_i\Delta t)}\bigg)^2.
\end{align*}
We will next demonstrate how to find explicit formulas for the partial derivatives of the distance functions. There exists a rotation of angle $\theta$, denote as $R_{\theta}$, such that \[R_\theta(\Omega_i(\sigma_i\mu_i\Delta t))=\{(x,y):y\geq f(x) \} \] for $f(x)=\sqrt{\sigma_i\mu_i\Delta t}\phi(|x|/\sqrt{\sigma_i\mu_i \Delta t})$. For the moment, consider the case where $\theta=0$. Let 
\[
g(x,y,p)=\frac{-(x-p)f'(p)+(y-f(p))}{\sqrt{1+(f'(p))^2}}.
\]
Additionally, let $p^*(x,y)=\displaystyle \argmin_{p} (x-p)^2+(y-f(p))^2$ be the $x$-coordinate of the closest point on the curve $(q,f(q))$ to $(x,y)$. Note that \begin{align}
\label{eq:der}
\frac{d}{dp}[(x-p)^2+(y-f(p))^2]|_{p^*} =2(x-p^*)+2(y-f(p^*))f'(p^*)=0.\end{align}
 We have $g(x,y,p^*(x,y))=d_{\Omega_i(\sigma_i\mu_i\Delta t)}(x,y)$ (see proposition 1 in~\cite{esedog2010diffusion}). We can find the partial derivatives of the distance function by differentiating $g$, for example: 
\begin{align*}
&\frac{d}{dx} d_{\Omega_i(\sigma_i\mu_i\Delta t)}=g_x+g_pp^*_x\\
&\frac{d^2}{dx^2} d_{\Omega_i(\sigma_i\mu_i\Delta t)}=g_{xx}+2g_{px}p^*_x+g_{pp}(p^*_x)^2+g_{p}p^*_{xx}.
\end{align*}
The partial derivatives of $p^*$ are obtained by implicitly differentiating  $(x-p^*)+(y-f(p^*))f'(p^*)=0$. We have that 
\begin{align*}
&g_p=\frac{[(x-p)+(y-f(p))f'(p)]f''(p)}{(1+(f'(p))^2)^{3/2}}=0
\end{align*}
by \cref{eq:der}, so we do not need to solve for $p^*_{xx}$. Applying the above for an arbitrary angle $\theta$, the partial derivatives of the distance function are 
\small
\begin{align*}
&\frac{\partial}{\partial x} d_{\Omega_i(\sigma_i\mu_i\Delta t)}=\frac{-\cos(\theta)f'(p^*)+\sin(\theta)}{\sqrt{1+(f'(p^*))^2}}\\
&\frac{\partial}{\partial y} d_{\Omega_i(\sigma_i\mu_i\Delta t)}=\frac{\sin(\theta)f'(p^*)+\cos(\theta)}{\sqrt{1+(f'(p^*))^2}}\\
&\frac{\partial^2}{\partial x^2} d_{\Omega_i(\sigma_i\mu_i\Delta t)}=[-\cos^2(\theta)-2\sin(\theta)\cos(\theta)f'(p^*)-\sin^2(\theta)(f'(p^*))^2]h(x,y,p^*)\\
&\frac{\partial^2}{\partial x \partial y} d_{\Omega_i(\sigma_i\mu_i\Delta t)}=[\sin(\theta)\cos(\theta)-\cos(2\theta)f'(p^*)-\sin(\theta)\cos(\theta)(f'(p^*))^2]h(x,y,p^*)\\
&\frac{\partial^2}{\partial y^2} d_{\Omega_i(\sigma_i\mu_i\Delta t)}=[-\sin^2(\theta)+2\sin(\theta)\cos(\theta)f'(p^*)-\cos^2(\theta)(f'(p^*))^2]h(x,y,p^*)\\
&h(x,y,p^*)=\bigg(\frac{1}{1+(f'(p^*))^2+[f(p^*)-(x\sin(\theta)+y\cos(\theta))]f''(p^*)}\bigg)\bigg(\frac{f''(p^*)}{[1+(f'(p^*))^2]^{3/2}}\bigg).
\end{align*}
\normalsize
We now have explicit formulas for every step of Newton's method allowing us to quickly and accurately minimize \cref{eq:dist3d} to find the closest point projection.

\subsection{The DMIIM's Relationship to the VIIM}
\label{DMVIIM} 
In this section, we discuss the precise relationship between the DMIIM and the VIIM. We consider the very special case of equal surface tensions, $\sigma_i$ = 1 (for all $i$), corresponding to the Herring angles of $(120^\circ, 120^\circ, 120^\circ)$. This is the one case in which our numerical results from previous sections suggest the VIIM converges to the correct solution. The standard maximum principle for two-phase motion by mean curvature implies that overlaps cannot occur at the end of the curve evolution step of the VIIM or the DMIIM. The only interesting question is how the two algorithms allocate points in the ``vacuum'' region,  $\{\mathbf{z}:d_{\Sigma_j(\Delta t)}(\mathbf{z}) < 0 \text{ for all } j\}$.\\

Due to the symmetry of this situation, if $\mathbf{p}\in S$, then any $\mathbf{q}\in\mathbb{R}^3$ obtained by a permutation of the components of $\mathbf{p}$ also satisfies $\mathbf{q}\in S$.\\

Let $\mathbf{x} \in\mathbb{R}^3$ be given, with $\mathbf{x}_i < 0$ for all $i$.
Let $\mathbf{p} = \Pi_S(\mathbf{x})$ with $\mathbf{p}=\Phi_S(\mathbf{z})$ for some $\mathbf{z}\in\mathbb{R}^2$.
We are thus assuming implicitly that $\Pi_S(\mathbf{x})$ consists of a single point $\mathbf{p} \in S$.

\begin{claim}
\label{claim:lemma}
$\mathbf{x}_i =\max(\mathbf{x}_1,\mathbf{x}_2,\mathbf{x}_3)$ if and only if $\mathbf{p}_i=\max(\mathbf{p}_1,\mathbf{p}_2,\mathbf{p}_3)$

\end{claim}
\begin{proof}
The proof will be broken up into three parts: For $i\neq j$
(1) if $\mathbf{x}_i \geq \mathbf{x}_j$ then $\mathbf{p}_i \geq \mathbf{p}_j$,
(2) if $\mathbf{p}_i=\mathbf{p}_j$, then $\mathbf{x}_i=\mathbf{x}_j$ and
(3) if $\mathbf{p}_i \geq \mathbf{p}_j$, then $\mathbf{x}_i \geq \mathbf{x}_j$. Statements (1) and (3) then imply the claim.\\

(1) If $\mathbf{x}_i \geq \mathbf{x}_j$ and $\mathbf{p}_i < \mathbf{p}_j$, then
\[ ||\mathbf{x} - \mathbf{q}|| \leq ||\mathbf{x} - \mathbf{p}|| \]
where $\mathbf{q}\in S$ is obtained from $\mathbf{p}$ by interchanging its $i$-th and $j$-th components (since $\mathbf{p}_i \neq \mathbf{p}_j$, then $\mathbf{q}\neq \mathbf{p}$).
Indeed,
\[ ||\mathbf{x}-\mathbf{p}||^2 = ||\mathbf{x}-\mathbf{q}||^2 + 2(\mathbf{x}_i-\mathbf{x}_j)(\mathbf{p}_j-\mathbf{p}_i) \geq ||\mathbf{x}-\mathbf{q}||^2.\]
This contradicts $\mathbf{p} = \Pi_S(\mathbf{x})$; so the first statement is established.\\

(2) If $\mathbf{p}_i=\mathbf{p}_j$, let $\mathbf{n}$ denote a unit normal to $S$ at $\mathbf{p}$.
A short calculation shows that $\mathbf{n}_i = \mathbf{n}_j$.
Indeed: \[\mathbf{n}_i = \mathbf{n}_j=D_{\mathbf{u}}d_{\Omega_i(\Delta t)}|_{\Phi^{-1}(\mathbf{p})}\] where $\mathbf{u}$ is the unit vector perpendicular to $\partial \Omega_i(0) \cap \partial \Omega_j(0)$ pointing into $\Omega_i(0)$.
Since $\mathbf{p} = \Pi(\mathbf{x})$ implies that $\mathbf{x}-\mathbf{p} = \beta \mathbf{n}$ for some $\beta\in\mathbb{R}$, we get $\mathbf{x}_i = \mathbf{x}_j$.\\

(3) Assume $\mathbf{x}_i<\mathbf{x}_j$. Since $\mathbf{x}_i\leq\mathbf{x}_j$, by statement (1) $\mathbf{p}_i \leq \mathbf{p}_j$. Furthermore, since  $\mathbf{x}_i \neq \mathbf{x}_j$ statement (2) implies $\mathbf{p}_i \neq \mathbf{p}_j$. Hence $\mathbf{p}_i<\mathbf{p}_j$.
\end{proof}

\begin{claim}
Let the phases $\Sigma_i^k$ at time step $k$ have smooth boundaries, with triple junctions in the same neighborhood of $(120^\circ , 120^\circ , 120^\circ)$ as in \cref{claim:cond}.
Then, for every small enough time step size $\Delta t>0$, the DMIIM and the VIIM yield the same $\Sigma_i^{k+1}$.
\end{claim}
\begin{proof}
Due to the symmetry $\mathbf{z}\in \Omega_i(0)$ if and only if $d_{\Omega_i(\Delta t)}(\mathbf{z})=\max_j{d_{\Omega_j(\Delta t)}(\mathbf{z})}$. A consequence is that \begin{align}
\label{eq:equal}
\Phi(\Omega_i(0))=\{\mathbf{p}\in S:\mathbf{p}_i=\max(\mathbf{p})\}.
\end{align}
The projection map $\Phi$ is well defined for points $(d_{\Sigma_1(\Delta t)}(\mathbf{z}),d_{\Sigma_2(\Delta t)}(\mathbf{z}),d_{\Sigma_3(\Delta t)}(\mathbf{z}))$, where $d_{\Sigma_j(\Delta t)}(\mathbf{z})< 0$ for all $j$ by \cref{claim:cond}. Then \cref{claim:lemma} along with \cref{eq:equal} give us that $d_{\Sigma_i(\Delta t)}(\mathbf{z})=\max_j{d_{\Sigma_j(\Delta t)}(\mathbf{z})}$ if and only if
\[
\Pi_S(d_{\Sigma_1(\Delta t)}(\mathbf{z}),d_{\Sigma_2(\Delta t)}(\mathbf{z}),d_{\Sigma_3(\Delta t)}(\mathbf{z}))\in \Phi(\Omega_i(0))
\]
completing the proof.
\end{proof}

The dictionary mapping implicit interface method is then an extension of the Voronoi implicit interface method to cases of unequal surface tension. As we show in the next section, the DMIIM numerically converges to the exact solution in cases of unequal surface tensions.

\subsection{Numerical Results for the DMIIM}
\begin{algorithm}
  \caption{Parameterized DMIIM for ``Grim Reaper'' tests with $\theta_2=\theta_3$.
    \label{alg:DMIIMpar}}
  \begin{algorithmic}[1]
  \STATE Let $N=T/\Delta t$.
    \STATE Choose points $\{x_i\}_{i=1}^n \in [0,.25]$ and set $y_i^0=f(x_i,0)$.
      \FOR{$k \gets 1 \textrm{ to } N$}
          \STATE Use $\{x_i,y_i^{k-1}\}_{i=1}^n$ to parameterize $\partial\Sigma_1$ and $\partial\Sigma_2$, denoted as $\gamma^{+}$ and $\gamma^{-}$ respectively.
          \STATE Evolve $\gamma^{+}$ and $\gamma^{-}$ by \cref{eq:curve} for time $\mu_1\sigma_1\Delta t$ and $\mu_2\sigma_2\Delta t$ respectively.
          \STATE For each $x_i$ find $\tilde{y}_i$ such that
    \begin{align*}
    (x^*,y^*)=\Phi^{-1} \circ \Pi_{S}(d_{\Sigma_1(\sigma_1\mu_1\Delta t)}(x_i,\tilde{y}_i),d_{\Sigma_2(\sigma_2\mu_2\Delta t)}(x_i,\tilde{y}_i),d_{\Sigma_3(\sigma_3\mu_3\Delta t)}(x_i,\tilde{y}_i))
    \end{align*}  satisfies $(x^*,y^*) \in \partial \Omega_1(0) \cup \partial \Omega_2(0)$.
    \STATE $y^{k}_i \gets \tilde{y}_i$
      \ENDFOR
  \end{algorithmic}{}
\end{algorithm}

\label{sec:nrdm}

\begin{table}
\parbox{.49\linewidth}{
\begin{center}
\caption{DMIIM\\
$\theta_1=120^\circ$}
\label{tab:dict120}
\begin{tabular}{|c|c|c|c|}
\hline
$\Delta t$&$n$&RE&Order\\
\hline
$2^{-10}$&$1024$&0.0202&-\\
\hline
$2^{-11}$&$1449$&0.0141&0.520\\
\hline
$2^{-12}$&$2048$&0.0099&0.513\\
\hline
$2^{-13}$&$2897$&0.0069&0.509\\
\hline
$2^{-14}$&$4096$&0.0049& 0.507\\
\hline
$2^{-15}$&$5793$&0.0034&0.505\\
\hline
$2^{-16}$&$8192$&0.0024&0.503\\
\hline
$2^{-17}$&$11586$&0.0017&0.498\\

\hline
\end{tabular}
\end{center}
}
\parbox{.49\linewidth}{
\begin{center}
\caption{DMIIM\\ $\theta_1=90^\circ$}
\label{tab:dict90}
\begin{tabular}{|c|c|c|c|}
\hline
$\Delta t$&$n$&RE &Order\\
\hline
$2^{-10}$&$1024$&0.00207&-\\
\hline
$2^{-11}$&$1449$&0.00107&0.954\\
\hline
$2^{-12}$&$2048$&0.00056&0.922\\
\hline
$2^{-13}$&$2897$&0.00031&0.842\\
\hline
$2^{-14}$&$4096$&0.00018&0.772\\
\hline
$2^{-15}$&$5793$&0.00011&0.701\\
\hline
$2^{-16}$&$8192$&0.00007&0.614\\
\hline
$2^{-17}$&$11586$&0.00005&0.576\\
\hline
\end{tabular}
\end{center}
}
\parbox{.49\linewidth}{
\begin{center}
\captionsetup{justification=centering}
\caption{DMIIM\\
$(\theta_1,\theta_2,\theta_3)=(75^\circ,135^\circ,150^\circ)$ with $\mu_i=1$}
\label{tab:DMIIM75mu1}
\begin{tabular}{|c|c|c|c|}
\hline
$\Delta t$&$n$&RE &Order\\
\hline
$2^{-10}$&$1024$&0.0067&-\\
\hline
$2^{-11}$&$1449$&0.0053&0.338\\
\hline
$2^{-12}$&$2048$&0.0040&0.411\\
\hline
$2^{-13}$&$2897$&0.0029&0.450\\
\hline
$2^{-14}$&$4096$&0.0021&0.470\\
\hline
$2^{-15}$&$5793$&0.0015&0.480\\
\hline
$2^{-16}$&$8192$&0.0011&0.494\\
\hline
\end{tabular}
\end{center}
}
\parbox{.49\linewidth}{
\begin{center}
\captionsetup{justification=centering}
\caption{DMIIM\\ 
$(\theta_1,\theta_2,\theta_3)=(75^\circ,135^\circ,150^\circ)$ with $\mu_i=\frac{1}{\sigma_i}$}
\label{tab:DMIIM75musig}
\begin{tabular}{|c|c|c|c|}
\hline
$\Delta t$&$n$&RE &Order\\
\hline
$2^{-10}$&$1024$&0.0138&-\\
\hline
$2^{-11}$&$1449$&0.0094&0.548\\
\hline
$2^{-12}$&$2048$&0.0065&0.540\\
\hline
$2^{-13}$&$2897$&0.0045&0.533\\
\hline
$2^{-14}$&$4096$&0.0031&0.527\\
\hline
$2^{-15}$&$5793$&0.0022&0.522\\
\hline
$2^{-16}$&$8192$&0.0015&0.521\\
\hline
\end{tabular}
\end{center}
}
\end{table}
In this section, we perform on the DMIIM the same careful numerical convergence tests that we subjected the VIIM to. In addition we test on some examples where all the surface tensions are different.
\cref{alg:DMIIMpar} details the implementation of the DMIIM with parameterized curves for ``Grim Reaper''  tests with $\theta_2=\theta_3$. The implementation of the non-symmetric case is similar. We note that in our implementation of \cref{alg:DMIIMpar} that projecting onto the template surface is the computation bottleneck, taking about 25 times longer than the Voronoi reconstruction step it replaces in the VIIM. This is mainly due to the number of points we use to represent each $\Omega_i$. We choose to err on the side of caution in our numerical studies, representing $\Omega_i$ with many more points than needed.\\

 We run the following ``Grim Reaper'' tests:
\begin{enumerate}
\item Angles $(\theta_1,\theta_2,\theta_3)=(120^\circ,120^\circ,120^\circ)$ with $\sigma_1=\sigma_2=\sigma_3=1$.
\item Angles $(\theta_1,\theta_2,\theta_3)=(90^\circ,135^\circ,135^\circ)$ with $\sigma_1=2-\sqrt{2}$ and $\sigma_2=\sigma_3=\sqrt{2}$
\item Angles $(\theta_1,\theta_2,\theta_3)=(75^\circ,135^\circ,150^\circ)$ with $\sigma_1=\frac{1}{4}(-2+3\sqrt{2}+\sqrt{6})$, $\sigma_2=\frac{1}{4}(2-\sqrt{2}+\sqrt{6})$, $\sigma_3=\frac{1}{4}(2+\sqrt{2}-\sqrt{6})$ and $\mu_i=1$.
\item Angles $(\theta_1,\theta_2,\theta_3)=(75^\circ,135^\circ,150^\circ)$  with the same $\sigma_i$'s as the previous test with $\mu_i=\frac{1}{\sigma_i}$.
\end{enumerate}
Each of the simulations use the following parameters:
\begin{itemize}
\item The total time the system is evolved: $T=18/512$.
\item Number of points tracked on the parameterized curve: $n$ as given in the table.
\item Step size in \cref{eq:curve}: $\delta t=\sigma\Delta t/n$.
\item We are measuring the relative error (RE) of the area of symmetric difference of phase $\Sigma_1$ in $\{(x,y):0\leq x\leq \beta-.04 \text{ or } \beta+.04\leq x\leq .5\}$, see \cref{rel}.
\end{itemize}
The results are contained in \cref{tab:dict120} through \cref{tab:DMIIM75musig}.
\subsection{Level Set Examples of the DMIIM}
\label{sec:four}
We demonstrate the level set formulation of the DMIIM in two and three dimensions. In 2d we evolve a system that goes through a well understood topological change. Initially, we have two ``Grim Reapers'' translating vertically towards each other until they collide. After the collision, two new junctions form that travel horizontally away from each other forming a new horizontal interface between the top and bottom phases.\\

\begin{figure}

  \begin{center}
\includegraphics[width=.24\textwidth]{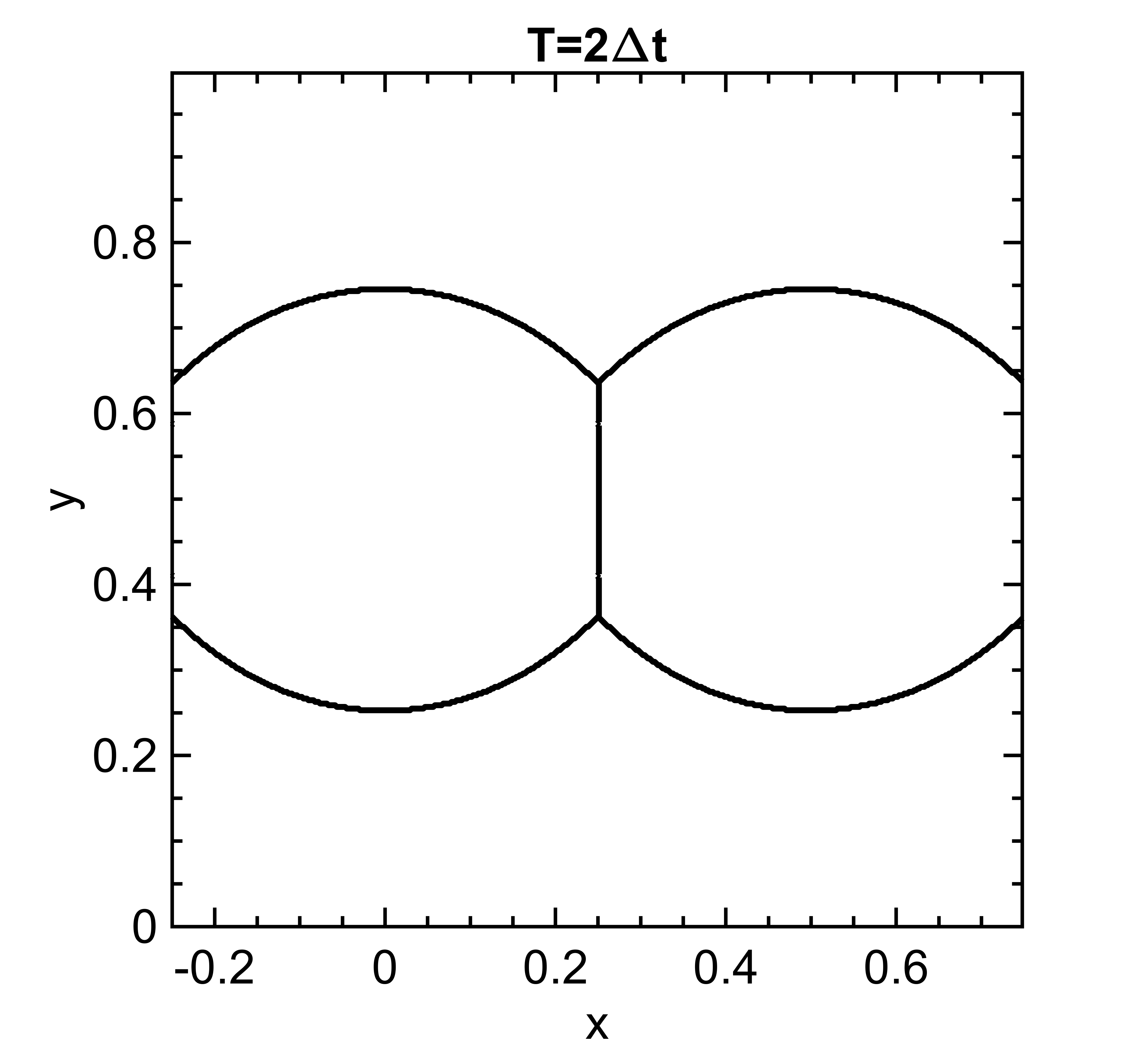} %
\includegraphics[width=.24\textwidth]{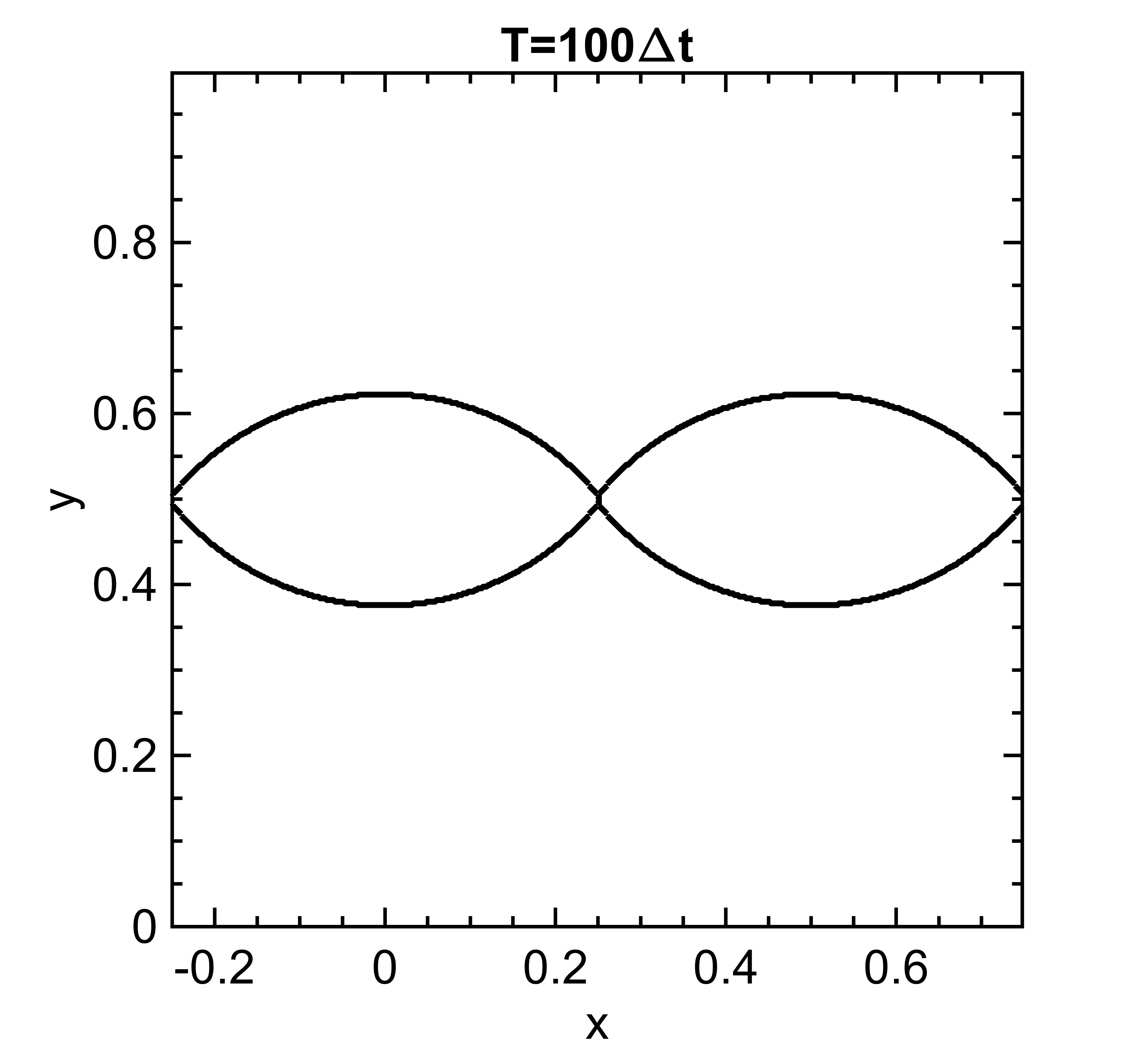}%
\includegraphics[width=.24\textwidth]{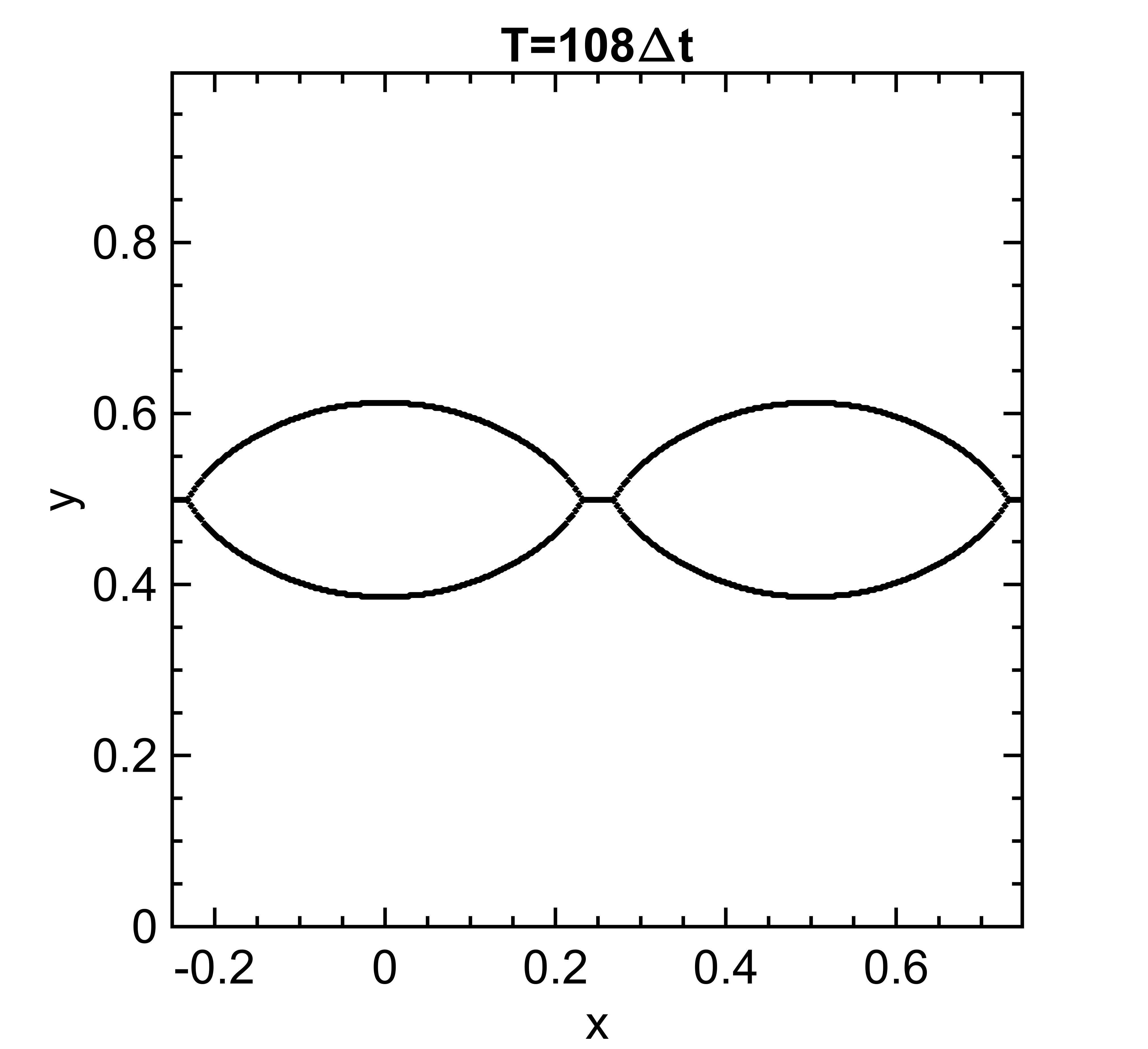}%
\includegraphics[width=.24\textwidth]{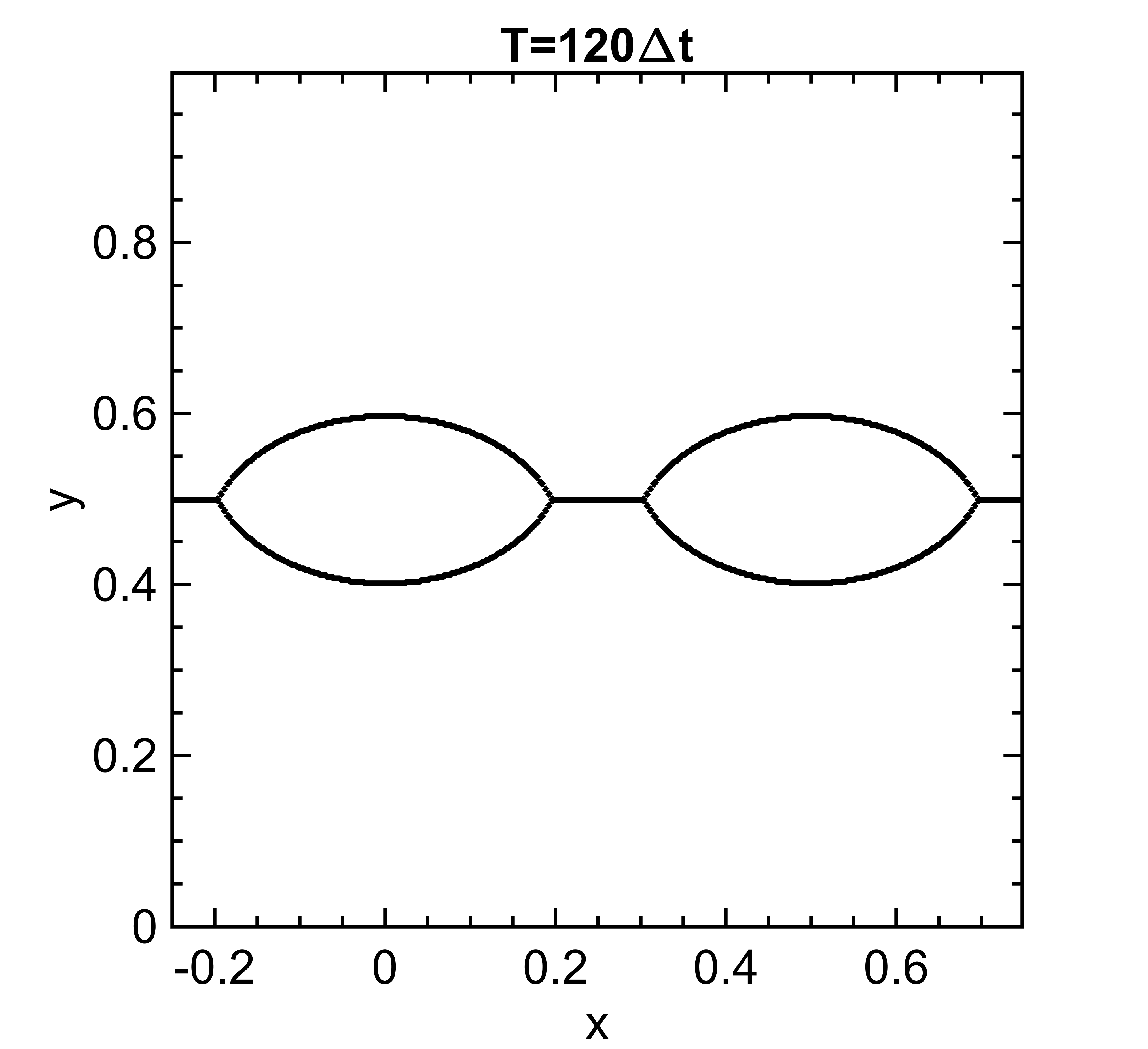}%
\end{center}
\caption{\footnotesize Two ``Grim Reapers'' colliding, computed using the DMIIM algorithm (practical implementation on uniform grid). The initial $(90^\circ,135^\circ,135^\circ)$ angles change to $(120^\circ,120^\circ,120^\circ)$ after the collision, as expected. The multiple junction that momentarily forms at collision appears to be handled appropriately.}
\label{fig:split}
\end{figure}

The initial configuration is 
\footnotesize
\begin{align*}
\Sigma^0_1=&\bigg\{(x,y): y>\frac{3}{4}+f(\frac{1}{4}-|x-\frac{1}{4}|)\bigg\}\\
\Sigma^0_2=&\bigg\{(x,y):x<\frac{1}{4} \text{ and } \frac{1}{4}+f(x)<y<\frac{3}{4}+f(x)\bigg\}\\
\Sigma^0_3=&\bigg\{(x,y):x<\frac{1}{4} \text{ and }\frac{1}{4}+f(x)<y<\frac{3}{4}+f(x)\bigg\}\\
\Sigma^0_4=&\bigg\{(x,y):\frac{1}{4}+f(\frac{1}{4}-|x-\frac{1}{4}|)\bigg\}
\end{align*}
\normalsize
where $f(x)=\frac{1}{\pi}\log(\cos(\pi x)).$
We use the surface tensions matrix
\footnotesize
\[
\sigma=\left(\begin{array}{cccc}
0 & 1 & 1 & 1\\
1 & 0 & \sqrt{2} & 1\\
1 & \sqrt{2} & 0 & 1\\
1 & 1 & 1 & 0
\end{array}\right)
\]
\normalsize
and $(\mu_1\sigma_1,\mu_2\sigma_2,\mu_3\sigma_3,\mu_4\sigma_4)=(2-\sqrt{2},\sqrt{2},\sqrt{2},2-\sqrt{2})$.
Before the topological change the angles are $(90^\circ,135^\circ,135^\circ)$ at the triple junctions. After the change the angles are $(120^\circ,120^\circ,120^\circ)$. For the simulation we use the domain $[0,1]\times[0,1]$ with periodic boundary conditions on a 512 by 512 grid. We set $\Delta t=\frac{1}{1600\pi\sin(3\pi/4)}$. \cref{fig:split} shows this system at different times.\\

\begin{figure}

  \begin{center}
\includegraphics[width=.49\textwidth]{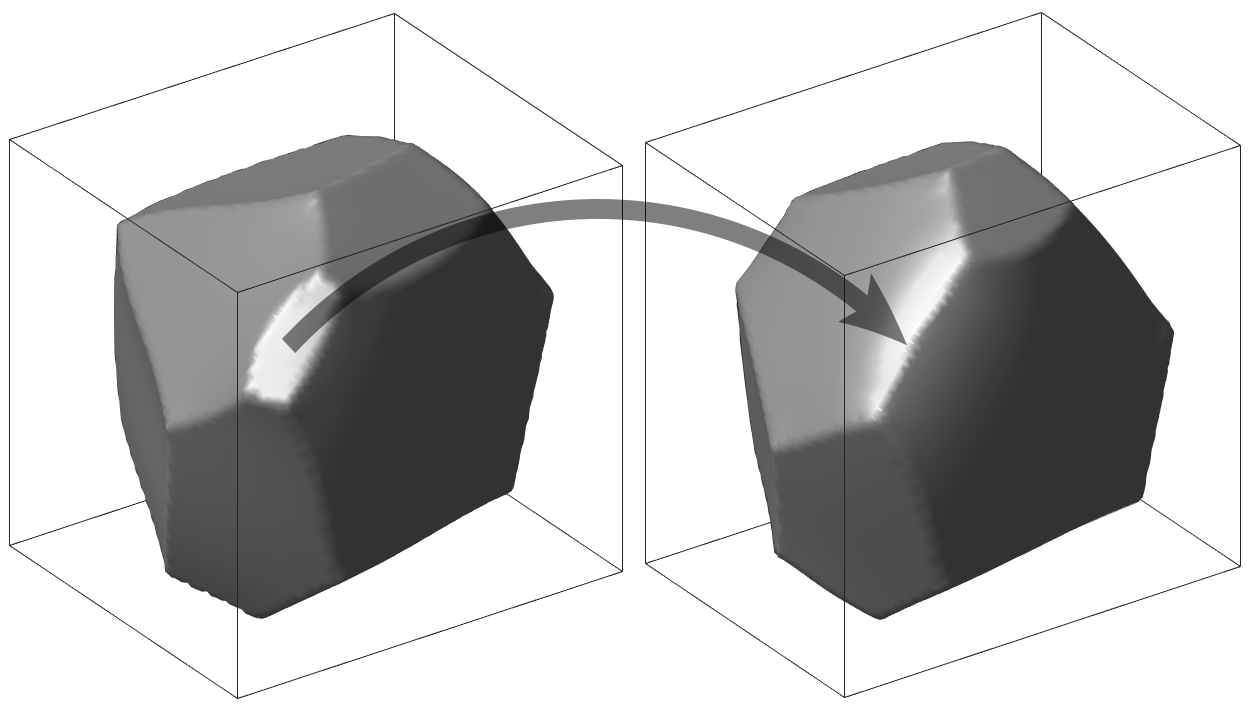} %
\includegraphics[width=.49\textwidth]{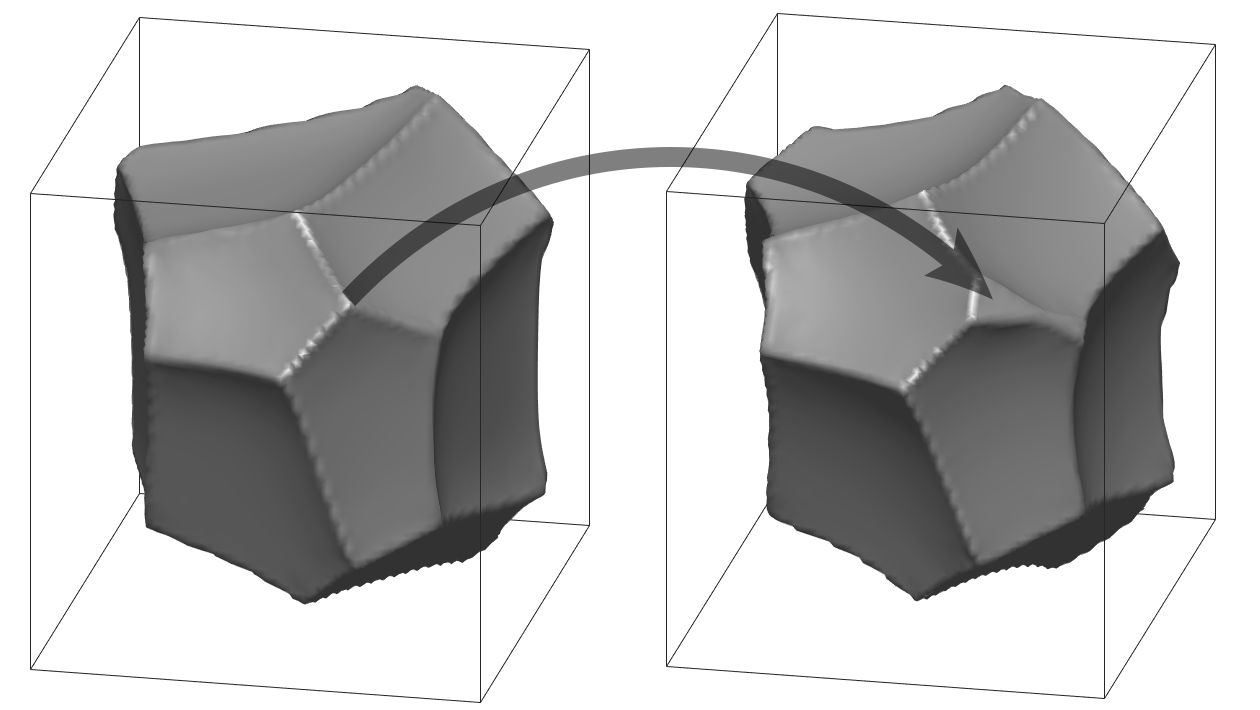}%
\end{center}
\caption{\footnotesize Two grains undergoing topological
changes on the face of the grains. In the first grain, four quadruple points collide. In the second grain, a quadruple point splits.}
\label{fig:3d}
\end{figure}

In three dimensions, we evolve a system starting from a Voronoi diagram of 8 points
taken at random on the 3-torus. Six of the phases have surface tension equal to $1$ while the other two phases have surface tension $\sqrt{2}-1$. Thus three angle configurations are possible: $(120^\circ,120^\circ,120^\circ)$, $(90^\circ,135^\circ,135^\circ)$, and $\approx (146^\circ,107^\circ,107^\circ)$ Quadruple points split and collide throughout the
evolution causing the faces of the grains to undergo topological changes as seen in \cref{fig:3d}. The two preceding examples show reasonable behavior and demonstrate the practical use of the DMIIM as a level set method.


\section{Conclusion}
\label{sec:conclusions}
In this work, we have presented careful numerical convergence studies showing that the Voronoi implicit interface method gets close but does not converge to the correct evolution in the unequal surface tension case of multiphase motion by mean curvature.
In addition, we proposed a correction to the method that fixes the non-convergence while maintaining the simplicity and the spirit of the method; indeed, the new algorithm reduces to the original in the case of equal surface tensions.
The correction is in the spirit of the projection method \cite{ruuth} of Ruuth proposed in the context of threshold dynamics.
We subjected the new algorithm to the same rigorous numerical convergence studies as the original, verifying that the non-convergence of the latter is rectified.
As in~\cite{ruuth}, the new algorithm is somewhat heuristic in the handling of higher order junctions ($\geq 4$) -- but numerical evidence is presented that suggests it behaves reasonably even in their presence.
Nevertheless, a more systematic approach, perhaps a variational interpretation of the VIIM in the spirit of \cite{esedog2015threshold}, would be far preferable, not least as a more reliable extension to arbitrary junctions.
Highlighting this need for further investigation of the VIIM is perhaps the most notable contribution of the present study.\\
\bibliographystyle{siamplain}
\bibliography{references}
\end{document}